\documentclass[12pt, a4paper]{article}

\usepackage{amsmath, amsthm, amssymb}  
\usepackage{fancyhdr}
\usepackage{graphicx}
\usepackage{enumerate}
\usepackage[top=2.5cm, bottom=2.5cm, left=2.5cm, right=2.5cm]{geometry}

\usepackage{mathrsfs} 
\usepackage{mathtools}
\usepackage{color}

\newtheorem*{definition}{Definition}{\rm}
\newtheorem*{definition*}{Definition} 
\newtheorem{lemma}{\bf Lemma}
\newtheorem{proposition}{\bf Proposition}
\newtheorem{theorem}{\bf Theorem}
\newtheorem{corollary}{\bf Corollary}

\renewenvironment{proof}{\noindent {\it Proof: }}{\rm\\}
\theoremstyle{definition}
\newtheorem*{remark}{\it Remark}{\rm}
\newtheorem*{example}{\it Example}{\rm}

\usepackage{algorithm}
\usepackage{algorithmic}
\usepackage[T1]{fontenc}

\makeatletter
\newsavebox\myboxA
\newsavebox\myboxB
\newlength\mylenA

\newcommand*\xoverline[2][0.75]{%
    \sbox{\myboxA}{$\m@th#2$}%
    \setbox\myboxB\null
    \ht\myboxB=\ht\myboxA%
    \dp\myboxB=\dp\myboxA%
    \wd\myboxB=#1\wd\myboxA
    \sbox\myboxB{$\m@th\overline{\copy\myboxB}$}
    \setlength\mylenA{\the\wd\myboxA}
    \addtolength\mylenA{-\the\wd\myboxB}%
    \ifdim\wd\myboxB<\wd\myboxA%
       \rlap{\hskip 0.5\mylenA\usebox\myboxB}{\usebox\myboxA}%
    \else
        \hskip -0.5\mylenA\rlap{\usebox\myboxA}{\hskip 0.5\mylenA\usebox\myboxB}%
    \fi}
\makeatother

\DeclareFontFamily{U}{mathx}{\hyphenchar\font45}
\DeclareFontShape{U}{mathx}{m}{n}{ <-> mathx10 }{}
\DeclareSymbolFont{mathx}{U}{mathx}{m}{n}
\DeclareFontSubstitution{U}{mathx}{m}{n}
\DeclareMathAccent{\widebar}{\mathalpha}{mathx}{"73}

\makeatletter
\newcommand{\cwidebar}[2][0]{{\mathpalette\@cwidebar{{#1}{#2}}}}
\newcommand{\@cwidebar}[2]{\@cwideb@r{#1}#2}
\newcommand{\@cwideb@r}[3]{%
  \sbox\z@{$\m@th#1\mkern-#2mu#3\mkern#2mu$}%
  \widebar{\box\z@}%
}
\makeatother

\usepackage{auto-pst-pdf}
\usepackage{pst-plot}
\usepackage{pdftricks}
\begin{psinputs}
\usepackage{pstricks,pst-plot}
\end{psinputs}

\newcommand{\closure}[2][5]{%
{}\mkern#1mu\overline{\mkern-#1mu#2}}

\begin{document}

\title{Topological spaces satisfying a closed graph theorem}
\author{Dominikus Noll
\thanks{Universit\'e de Toulouse, Institut de Math\'ematiques,  
118, route de Narbonne, 31062 Toulouse, France \newline
{\tt $\quad$ dominikus.noll$@$math.univ-toulouse.fr}}
}

\date{}

\maketitle

\begin{abstract}
We discuss topological versions of the closed graph theorem, where continuity is inferred from near continuity in tandem with suitable conditions on source or target spaces. 
We seek internal characterizations of spaces  satisfying
a closed graph theorem, and we compare closed graph and open mapping spaces.

\vspace{.2cm}
\centerline{\bf Key Words}
\noindent
Closed graph theorem $\cdot$  Open mapping theorem $\cdot$ Michael game $\cdot$ Hereditary Baire space $\cdot$ Banach-Mazur game $\cdot$ Tandem game

\vspace{.2cm}
\centerline{\bf MSC 2020}

\centerline
{54E52 $\cdot$  54A10  $\cdot$  54C08 $\cdot$ 54C10}
\end{abstract}

\section{Introduction}

By the term
{\it closed graph theorem} we understand any instance of the  following general
scheme:

{\it  Let $f$ be a nearly continuous mapping with closed graph  from a space $F$ to a space $E$. Then, under suitable conditions on $F$ and $E$, the mapping $f$ is continuous.} 

Pettis \cite{pettis} proved a first purely topological  closed graph theorem where $F,E$ are completely metrizable, expanding  on Banach's classical antecedent.
The author extended this to the case where $F$ is Baire and $E$ is sieve-complete \cite{graph_theorem}, and Moors gave a further extension,
where $E$ is partition complete \cite{moors1}.  It is folklore that the closed graph theorem holds when $E$ is locally compact, and Mahavier's theorem extends this  to the case where 
$E$ is locally countably compact and $F$ is a Fr\'echet space.

In the classical categories of locally convex vector spaces, linear topological spaces, or topological groups, the closed graph theorem is often deduced
from the open mapping theorem. This is the case in Banach-Schauder theory \cite{koethe},  in V. Pt\'ak's approach to the closed graph theorem \cite{ptak,ptak_early},
and similarly so in the framework of topological vector spaces \cite{adasch}, and for topological groups \cite{husain}. While the methods used to realize this reduction may vary,
it is clear that both mapping theorems are closely related. By an
open mapping theorem we understand a statement of
the following
general form:

{\it Let $f:E \to F$ be a continuous nearly open bijection. Then, under suitable conditions on the spaces $E$ and $F$, the mapping $f$ is open.} 

The first purely topological open mapping theorem was proved by Weston \cite{We} and assumes a completely metrizable $E$,  again
expanding on Banach's classical result.  
Weston's  result was subsequently extended to \v{C}ech complete spaces   \cite{BP}, to
monotonically \v Cech complete $E$ in \cite{wilhelm_criteria}, to almost \v Cech-complete spaces $E$ in  \cite{on_the_theory,baire_category}, and recently 
to tandem-Baire spaces  \cite{noll_tandem}.
The following definition due to Husain \cite{husain}  has been useful in the study of open mapping theorems. 

\begin{definition}
{\bf (Open mapping space)}.
Let $\mathscr K$ be a class
of Hausdorff spaces. A Hausdorff space $E$ is called a $B_r(\mathscr K)$-space if the open mapping theorem holds between $E$ and  every $F\in \mathscr K$.
When $\mathscr K$ is the class of all
Hausdorff spaces, then $E$ is called a $B_r$-space. 
\end{definition}

The name is with reference to V. Pt\'ak's  open mapping theorem for  $B_r$-complete spaces \cite{ptak,ptak_early,koethe}.
The notion was also applied in topological vector spaces  \cite{adasch}, and topological groups  \cite{husain}. 
We take this as  motivation to propose the following

\begin{definition}
{\bf (Closed graph space)}.
Let $\mathscr K$ be a class of Hausdorff spaces. A Hausdorff space $E$ is called a $C(\mathscr K)$-space if the closed
graph theorem holds between every $F \in \mathscr K$ and $E$. We say $C$-space when $\mathscr K$ is the class of all Hausdorff spaces. 
\end{definition}

A variant of the closed graph theorem appears
in the work of Byzkowski and Pol \cite{BP}, Wilhelm \cite{wilhelm_relations,wilhelm_pettis} and Piotrowski and Szyma\'nski \cite{piotrowski}. These authors observed that 
closed graph group homomorphisms  exhibit  a somewhat stronger form of closedness, referred to as {\it separatingness}. Graph theorems based on separating functions $f$ 
therefore still generalize classical closed graph theorems, but lead to a richer theory in the topological context. Spaces $E$ satisfying the separated graph theorem will therefore be
termed $C_s(\mathscr K)$-spaces, the exact definition to be given in the next section. Properties of both classes $C(\mathscr K)$ and $C_s(\mathscr K)$
will be studied along with those of the closely related  $B_r(\mathscr K)$-spaces.

The structure of the paper is as follows. Section \ref{prelim} concretizes the set-up and gives preliminary results. In
Section \ref{enlarge} we show that a $C(\mathscr K)$-space remains $C(\mathscr K')$ for certain enlargements $\mathscr K'$ of $\mathscr K$. 
Section
\ref{sect_perfect} 
discusses invariance of $C(\mathscr K)$- and $C_s(\mathscr K)$-spaces under perfect preimages. Section \ref{sect_F}  inspects $F$-heredity,
and Section \ref{sect_product} examines product invariance.  Invariance properties  are used to
demarcate between  the classes $C(\mathscr K)$,
$C_s(\mathscr K)$ and $B_r(\mathscr K)$.
Section \ref{sect_M} studies
$C$- and $C(\mathscr M)$-spaces, where $\mathscr M$ means metrizable. 
Section \ref{sect_regular} proves that first countable $C_s$-spaces and first countable spaces $E \in \mathscr B \cap C_s(\mathscr B)$, $\mathscr B$ the Baire spaces in the sense of category,  
are regular, an observation which parallels the known fact
that $B_r$-spaces are semi-regular. In Section \ref{sect_tandem} we introduce a tandem version of the Michael game \cite{michael_game}, which leads to a new class of
hereditary Baire spaces, giving rise to an extension  of Moors' closed graph theorem \cite{moors1}. 
In Section \ref{sect_CM}  these spaces are then shown to belong to the class $C_s(\mathscr M_c)$, which leads to yet another  extension of
Pettis' closed graph theorem. The final Section \ref{sect_final} presents an internal characterization of the class of $C(\mathscr M_c)$-spaces,
where $\mathscr M_c$ stands for completely metrizable. This is motivated by a result of Kenderov and Revalski \cite{kenderov}, where the authors
characterize spaces $E$ containing a dense completely metrizable subspace by a variant of the closed graph theorem for set-valued mappings.

\section{Preliminaries}
\label{prelim}
All topological spaces considered are Hausdorff  and are often just referred to as {\it spaces}. 
A mapping $f:F \to E$  is {\it nearly continuous}  if $f^{-1}(V) \subseteq \xoverline[.9]{f^{-1}(V)}^\circ$ for every open $V\subseteq E$. Similarly,
$f:E \to F$ is {\it nearly open}  if $f(U) \subseteq \xoverline[.9]{f(U)}^\circ$ for every open $U\subseteq E$. 
General terminology follows \cite{engelking}.

We start our investigation with  
the following useful observation, which is  clear as continuous mappings between Hausdorff spaces have closed graph:

\begin{proposition}
\label{prop2}
Every $C(\mathscr K)$-space is a  $B_r(\mathscr K)$-space. 
\hfill $\square$
\end{proposition}
In other words, the open mapping theorem follows from the closed graph theorem. 
In the classical categories the converse is also true, see  Pt\'ak \cite{ptak}, \cite[p.29]{koethe} for locally convex vector spaces,  \cite{adasch} for topological vector spaces, 
and   \cite{husain} for topological groups.  We sketch a proof in the latter case:

\begin{proposition}
\label{groups}
Every $B_r$-group is a $C$-group.
\end{proposition}

\begin{proof}
Let $f:F \to E$ be a closed graph nearly continuous surjective homomorphism from the Hausdorff topological group $F$ onto the
$B_r$-group $E$. Let $\tau$ be the group topology on $E$ and define a coarser group topology $\tau'$ on $E$
by taking as a neighborhood base of the unit $e\in E$ all sets of the form
$V' = f(\xoverline[.9]{f^{-1}(V)}^\circ)^\circ$, where $V$ runs through the $\tau$-neighborhoods $V$ of $e$.  On checks that 
$f:F \to (E,\tau')$ is continuous,  that $id_E:(E,\tau) \to (E,\tau')$ is continuous and nearly open, and last but not least, that $\tau'$ is Hausdorff.
Then $i_E$ is open because $E$ is a $B_r$-group, so $\tau = \tau'$, hence $f:F \to (E,\tau)$ is continuous.
\hfill $\square$
\end{proof}

The argument makes heavy use of the group structure,  and it is not surprising that
there is no direct analogue in the topological setting.
Deriving a closed graph theorem from the open mapping theorem, if possible at all, 
must therefore be attempted differently.
One approach which narrows the gap between $C(\mathscr K)$- and $B_r(\mathscr K)$-spaces   uses a stronger notion of graph closedness proposed in
\cite{BP}, \cite{wilhelm_relations} and \cite{piotrowski}.

\begin{definition}
A mapping $f:F \to E$ is called {\rm separating} if for any $x,y\in E$ with $x\not= y$, there exist neighborhoods $V_x$ of $x$ and $V_y$ of $y$ in $E$
such that $\xoverline[.9]{f^{-1}(V_x)} \cap f^{-1}(V_y) = \emptyset = f^{-1}(V_x) \cap \xoverline[.9]{f^{-1}(V_y)}$.
\end{definition}

The following facts are known (see \cite{piotrowski,wilhelm_criteria,wilhelm_relations}): (a) Every closed graph group homomorphism is separating,
(b) every separating map has closed graph,  (c) every continuous $f$ is separating, (d) if $f$ has closed graph,
preimages of compact sets are compact, and $E$ is a $k$-space, then $f$ is separating,  (e)  there exist closed graph mappings which are not separating.
Separatingness is therefore a strengthened form of graph closedness, and both coincide in the classical categories.
This justifies the following

\begin{definition}
{\bf (Separatingness spaces)}.
Let $\mathscr K$ be a class of spaces,  then $E$ is called a $C_s(\mathscr K)$-space if every
nearly continuous separating mapping $f$ from any $F \in \mathscr K$ to $E$ is continuous. 
When $\mathscr K$
is the class of all Hausdorff spaces, we  say $C_s$-space.
\end{definition}

Clearly $C(\mathscr K)$-spaces are $C_s(\mathscr K)$-spaces.
Since the condition imposed on the graph of $f$ is stronger, requirements on $F,E$ may be weaker, so that 
the $C_s(\mathscr K)$-concept  leads to a wider range of applications.  
We start our investigation with
the following extension of  Proposition \ref{prop2}:

\begin{proposition}
\label{prop4}
Every $C_s(\mathscr K)$-space is a $B_r(\mathscr K)$-space.
\end{proposition}

\begin{proof}
Let $F \in \mathscr K$, and
let $f:E \to F$ be a continuous nearly open bijection.  We have to show that $f$ is open.
We show that $f^{-1}:F \to E$ is nearly continuous and separating. Since $E$ is a $C_s(\mathscr K)$-space, this implies continuity of $f^{-1}$,
hence openness of $f$.

Near continuity of $f^{-1}$ is obvious. Let $x\not= x'$ in $E$, $y=f(x)$,   $y'=f(x')$, then $(x,y')\not \in G(f)$, and since $f$ as a continuous function has closed graph,
there exist
$x\in U$ open in $E$, $y'\in V'$ open in $F$ with $(U \times V' )\cap G(f)=\emptyset$, or what is the same, $f(U) \cap V' = \emptyset$. Hence also $\xoverline[.9]{f(U)} \cap V'=\emptyset$.
Using continuity at $x'$ find $x'\in U'$ open with $f(U')\subseteq V'$. Then $\xoverline[.9]{f(U)} \cap f(U')=\emptyset$, and then of course also $\xoverline[.9]{f(U)}^\circ \cap f(U')=\emptyset$.
This gives $\xoverline[.9]{f(U)}^\circ \cap \xoverline[.9]{f(U')}=\emptyset$, and since $f$ is nearly open, 
$f(U) \subseteq  \xoverline[.9]{f(U)}^\circ$, hence $f(U) \cap \xoverline[.9]{f(U')} = \emptyset$. Altogether these prove that $f^{-1}$ is separating.
\hfill $\square$
\end{proof}

Before diving further into the analysis of
$C(\mathscr K)$- and $C_s(\mathscr K)$-spaces, let us mention that their definition may seem arbitrary in so far as we could just as well have kept  
$F$  fixed and let  $E$ vary in a class $\mathscr K$.
We would  call $F$ a $c(\mathscr K)$-space if the closed graph theorem  holds between $F$ and
every space $E\in\mathscr K$,  and similarly for
$c_s(\mathscr K)$. 
It turns out that this concept is of less interest.

\begin{proposition}
The following conditions on  $F$ are equivalent:
{\rm (i)} $F$ is a $c$-space,
{\rm (ii)} $F$ is a $c_s$-space,
{\rm (iii)} every dense subset $D$ in $F$ is open.
\end{proposition}

\begin{proof}
Clearly (i) implies (ii).  We show that  (ii) implies (iii). Let $\tau$ be the topology on $F$ and let $D \subseteq F$  be dense. 
Let $E$ be the point-set $F$ endowed with the finer topology generated by the basis  $\tau \cup \{U \cap D: U \in \tau\}$. Then $id_F:F \to E$ is separating
and nearly continuous, the latter since $U\cap D$ is dense in $U$. Since $E$ is a $c_s$-space, $id_F$ is continuous, hence $D \in \tau$. We show that (iii) implies (i). 
Let $f:F\to E$ be nearly continuous, then
for every open $U$ in $E$, $D = f^{-1}(U) \cup E \setminus \xoverline[.9]{f^{-1}(U)}$ is dense in $F$, hence open, hence
$\xoverline[.9]{f^{-1}(U)}^\circ \cap D$ is open, but this set equals $f^{-1}(U)$.
\hfill $\square$
\end{proof}

In consequence $c(\mathscr K)$- and $c_s(\mathscr K)$-spaces are not expected to lead to a reasonable theory when it comes to seeking necessary conditions, so that
we shall not employ these terms in the sequel.

\section{Enlargement  by taking images}
\label{enlarge}
In this section we show that if  $E$ is a $C(\mathscr K)$-space, then it remains a $C(\mathscr K')$-space for the larger class $\mathscr K'$
of feebly open quotients of spaces from $\mathscr K$. Recall that a mapping $f:F \to E$ is feebly open if $f(U)^\circ \not=\emptyset$ for every nonempty open $U \subseteq F$.

\begin{lemma}
\label{feeble}
Let $f:F \to E$ be nearly continuous with closed graph, and let $g:G \to F$ be continuous and feebly open. Then $f \circ g$ is nearly continuous and has closed graph.
Similarly, when $f:F \to E$ is nearly continuous and separating, then $f\circ g$ is also nearly continuous and separating.
\end{lemma}

\begin{proof}
When $f$ has closed graph,  then by continuity of $g$ the graph of $f \circ g$ is also closed.
Similarly, when  $f$ is separating, then $f \circ g$ is separating.

Concerning near continuity,
let $U\subseteq E$ be open.
We have to prove $g^{-1}(f^{-1}(U)) \subseteq \xoverline[.9]{g^{-1}(f^{-1}(U))}^\circ$. 
Since $f^{-1}(U) \subseteq \xoverline[.9]{f^{-1}(U)}^\circ$ and $g$ is continuous, it will be sufficient to prove
$g^{-1}(\xoverline[.9]{f^{-1}(U)}^\circ) \subseteq   \xoverline[.9]{g^{-1}(f^{-1}(U))}$.

Let $x$ be in the left hand side, $V$ an open  neighborhood of $x$. Since $g$ is continuous and $g(x) \in \xoverline[.9]{f^{-1}(U)}^\circ$, we may assume that
$g(V) \subseteq \xoverline[.9]{f^{-1}(U)}^\circ$. Since $g$ is feebly open, $g(V)^\circ$ is nonempty, and therefore $g(V) \cap f^{-1}(U)\not=\emptyset$,
because $f^{-1}(U)$ is dense in $ \xoverline[.9]{f^{-1}(U)}^\circ$. This shows $V \cap g^{-1}(f^{-1}(U)) \not=\emptyset$.
\hfill $\square$
\end{proof}

Given a class $\mathscr K$, let $\mathscr K'$ be the class of all feebly open quotients of spaces from $\mathscr K$. Then we have the following invariance result:

\begin{proposition}
\label{feebly}
Every $C(\mathscr K)$-space $E$ is also a $C(\mathscr K')$-space, and every $C_s(\mathscr K)$-space is a $C_s(\mathscr K')$-space.  \end{proposition}

\begin{proof}
Let $F \in \mathscr K'$ and let $f:F \to E$ be nearly continuous with closed graph. By definition of $\mathscr K'$ there exists $G \in \mathscr K$ and 
a continuous feebly open quotient map $g:G \to F$. By the Lemma
$f\circ g$ is nearly continuous with closed graph, hence continuous since $E$ is a $C(\mathscr K)$-space. Now as $g$ is quotient, $f$ is continuous.
The $C_s(\mathscr K)$-case is similar.
\hfill $\square$
\end{proof}

\begin{remark}
When $g$  is a continuous feebly open surjection, while not necessarily quotient, we can still achieve something.
Namely,  $f\circ g$ is still continuous, so for $U$ open 
$g^{-1}(f^{-1}(U))$ is open, hence by feeble openness of $g$ the image $g(g^{-1}(f^{-1}(U))) = f^{-1}(U)$ has nonempty interior, and that gives what is called {\it feeble continuity} of $f$.
The analogous result for the $C_s(\mathscr K)$-case is also true.
\end{remark}

The following is an immediate consequence of Proposition \ref{feebly}:

\begin{corollary}
Let $E$ be a $C(\mathscr K)$-space, and let $\mathscr K^\circ$ be the class of open continuous images of spaces from
$\mathscr K$. Then $E$ is a $C(\mathscr K^\circ)$-space. The analogous result for the $C_s(\mathscr K)$-case is also true. 
\hfill $\square$
\end{corollary}

Yet another consequence of Proposition \ref{feebly} is obtained when we recall that every
closed continuous and irreducible surjection is feebly open. 
Following
Gruenhage \cite{gruenhage} we call a continuous closed surjection $g:G \to F$ {\it inductively irreducible} if there exists
a closed subspace $H \subseteq G$ such that the restriction $h=g|H:H \to F$ is closed and irreducible.

The  classical result of
La\v snev \cite[5.5.12]{engelking}, which motors this topic,  says that  any continuous closed surjection $f$ from a paracompact space $G$
onto  a Fr\'echet space $F$ is inductively irreducible. In 
\cite{gruenhage} this is extended to $G$ meta-Lindel\"of and $F$ containing either a dense set of isolated points,
or a dense set of points which are limits of points taken from  a fixed  countable set. As further shown, when no weakly inaccessible cardinal exists, it suffices even to assume
that
$F$ does not contain any open dense-in-itself $P$-space.

\begin{proposition}
Let $\mathscr K$ be $F$-hereditary,  $\mathscr K_{ic}$  the class of  images of
spaces from $\mathscr K$ under inductively irreducible  mappings. 
Then every $C(\mathscr K)$-space $E$ is also a $C(\mathscr K_{ic})$-space. The same statement holds for $C_s(\mathscr K)$.
\end{proposition}

\begin{proof}
Let $F \in \mathscr K_{ic}$ and let $f:F \to E$ be nearly continuous with closed graph. By the definition of $\mathscr K_{ic}$ there exists $G \in \mathscr K$
and a continuous closed surjection $g:G \to F$ allowing a restriction $h:=g|H$ to a closed subspace $H$ of $G$ such that $h:H \to F$ is a closed irreducible surjection. Now 
$h$ is feebly open \cite{almost_complete}, as witnessed by
$\emptyset \not=F\setminus h(H\setminus U)\subseteq h(U)$ for nonempty open $U$. Hence by Lemma \ref{feeble}, $f\circ h$ has closed graph and is nearly continuous. Since $\mathscr K$ is $F$-hereditary,
we have $H \in \mathscr K$,  and since
$E$ is a $C(\mathscr K)$-space, $f\circ h$ is continuous. Since $h$ is quotient,  this implies continuity of $f$.
The $C_s(\mathscr K)$-case proceeds analogously. 
\hfill $\square$
\end{proof}

Recall that a continuous mapping is perfect if it is closed and has compact fibers; \cite[p. 182]{engelking}.  Now we have

\begin{corollary}
Let $\mathscr K$ be an $F$-hereditary  class, $\mathscr K_p$ the class of perfect images of spaces from $\mathscr K$. Then every
$C(\mathscr K)$-space is also a $C(\mathscr K_p)$-space. Analogously  in the $C_s(\mathscr K)$-case.
\end{corollary}

\begin{proof}
Perfect maps are inductively  irreducible  due to
compactness of the fibers.
\hfill $\square$
\end{proof}

\section{Invariance under perfect preimages}
\label{sect_perfect}
In this section we start discussing invariance properties of $C(\mathscr K)$- and $C_s(\mathscr K)$-spaces.
Recall that $F$ is called a $k$-space
if  any $A\subseteq F$ which has $A \cap K$ closed for every compact $K \subseteq F$ is closed in $F$; cf. \cite[p. 152]{engelking}.

\begin{proposition}
\label{perfect}
Let $\mathscr K$ be a class of $k$-spaces. Then
perfect preimages of $C(\mathscr K)$-spaces are $C(\mathscr K)$-spaces, and
perfect preimages of $C_s(\mathscr K)$-spaces are $C_s(\mathscr K)$-spaces.
\end{proposition}

\begin{proof}
1)
Let $p:E \to E_0$ be a perfect surjection, and suppose $E_0$ is a $C(\mathscr K)$-space.  We have to prove
that $E$ is a $C(\mathscr K)$-space.
Let $F \in \mathscr K$ and $f:F\to E$ nearly continuous with closed graph. We consider the mapping $f_0 = p \circ f:F \to E_0$.
As a composition of a nearly continuous and a continuous mapping $f_0$ is nearly continuous. Let us show that $f_0 = p \circ f$ has closed graph. 

Observe that $G(p \circ f)$ is the image of $G(f)$ under the cartesian product mapping ${\rm id}_F\times p: F \times E \to F \times E_0$.
Since $p$ is perfect, ${\rm id}_F \times p$ is closed (cf. \cite[Thm. 3.7.14]{engelking}), hence the image $G(p \circ f)$ of the closed set $G(f)$
under ${\rm id}_F \times p$ is closed in $F \times E_0$.
Since $E_0$ is a $C(\mathscr K)$-space, the mapping $f_0$ is continuous.

2) Now consider the case where $E_0$ is a $C_s(\mathscr K)$-space, and suppose $f:F \to E$ is nearly continuous and separating.  
We have to show that $f_0=p \circ f$
is separating. 
Consider $x,y\in E_0$, $x\not= y$, then $K_x=p^{-1}(x)$,
$K_y=p^{-1}(y)$ are disjoint compact sets in $E$, $p$ being perfect. Hence by \cite[Lemma 3]{piotrowski} there exist open sets
$K_x\subseteq V_x$, $K_y\subseteq V_y$ in $E$ such that $f^{-1}(V_x) \cap \xoverline[.9]{f^{-1}(V_y)} = \emptyset = \xoverline[.9]{f^{-1}(V_x)} \cap f^{-1}(V_y)$. Since $p$
is closed, there exist open sets $x\in U_x$, $y\in U_y$ such that $p^{-1}(U_x) \subseteq V_x$, $p^{-1}(U_y) \subseteq V_y$. Then
clearly $U_x,U_y$ are separating for $f_0=p\circ f$. To see the existence of $U_x$, 
note that $x \not\in p(E\setminus V_x)$ and this set is closed, hence there exists $U_x$ open with $x\in U_x$
and $U_x \cap p(E\setminus V_x)=\emptyset$. Then $p^{-1}(U_x) \subseteq V_x$. Similarly for $U_y$.
This shows continuity of $f_0$  for the $C_s(\mathscr K)$-case. 

3) Now let $K \subseteq F$ be compact. Then  as $f$ has closed graph, $f(K)$ is closed in $E$. 
But $f_0 = p\circ f$ is continuous, hence $p(f(K))$ is compact in $E_0$. Since $p$ is perfect, $p^{-1}(p(f(K)))$ is compact in $E$. Since $f(K)$ is closed in $E$,
it is also closed in $p^{-1}(p(f(K))$, hence is compact. It follows that $f$ maps compact sets to compact sets and has closed graph. 
Since $F\in \mathscr K$ is a $k$-space, this implies continuity of $f$.
\hfill $\square$
\end{proof}

An immediate consequence of Proposition \ref{perfect} is that if $\mathscr K$ is a class of $k$-spaces, $E$ is a $C(\mathscr K)$-space and $K$ is compact, then $E \times K$
is a $C(\mathscr K)$-space, and similarly for the $C_s(\mathscr K)$-case. 
Yet another consequence  is the following:

\begin{proposition}
\label{oplus}
Suppose $\mathscr K$ is a class of $k$-spaces and $E$ is a $C(\mathscr K)$-space,  then $E \oplus E$ is also a $C(\mathscr K)$-space.  The analogous result holds in the $C_s(\mathscr K)$-case.
\end{proposition}

\begin{proof}
Let $E \oplus E$ be realized as $(E \times \{1\}) \cup (E \times \{2\})$ and define $p:E \oplus E \to E$ as
$p(x,i) = x$. Then $p$ is perfect, hence the result follows with Proposition \ref{perfect}.
\hfill $\square$ 
\end{proof}

When $E$ is a $B_r$-space, then $E \oplus E$ is again $B_r$ even without the $k$-space restriction; cf. \cite{sums_and_products}. 
We do not know whether the $k$-space hypothesis can be avoided in Propositions \ref{perfect} and  \ref{oplus}.
Concerning sums we have the following

\begin{lemma}
\label{iota}
Let $(E_\iota: \iota \in I)$ with $|I|\geq 2$ be a family of spaces such that $E_\iota \oplus E_\kappa$ is a $C$-space for any two $\iota,\kappa \in I$, $\iota \not=\kappa$. Then
$E:=\bigoplus_{\iota \in I} E_\iota$ is a $C$-space. The analogous result holds for $C_s$-spaces. 
\end{lemma}

\begin{proof}
Let $f:F \to E = \bigoplus_\iota E_\iota$ be nearly continuous with closed graph. Let $F_\iota = f^{-1}(E_\iota)$, then
each restriction $f_\iota = f|F_\iota$, $f_\iota: F_\iota \to E_\iota$ has closed graph and is nearly continuous. Since as a closed-and-open subspace of the $C$-space $E_\iota \oplus E_\kappa$,  $E_\iota$ is 
also a $C$-space, the
$f_\iota$ are continuous.  Now for continuity of $f$ it suffices to prove that every $F_\iota$ is closed in $F$. Suppose this is not the case, then $\xoverline{F}_\iota \cap F_\kappa \not=
\emptyset$ for some $\kappa\not=\iota$, but that contradicts the hypothesis that $E_\iota \oplus E_\kappa$ is a $C$-space.  The proof of the $C_s$-case is analogous.
\hfill $\square$
\end{proof}

 By combining Corollary \ref{oplus} and the argument of Lemma \ref{iota}, we obtain

\begin{corollary}
\label{cor1} 
Let $\mathscr K$ be a hereditary class of $k$-spaces. Suppose $E$ is a $C(\mathscr K)$-space, then an arbitrary sum of copies of $E$ is a $C(\mathscr K)$-space.
The analogous result holds in the $C_s(\mathscr K)$-case.
\end{corollary}

Concerning sums we also have the following

\begin{proposition}
\label{prop10}
Let $E$ be a $C(\mathscr K)$-space, $L$ locally compact.  Suppose $\mathscr K$ is invariant under taking closed-and-open subsets. Then $E \oplus L$ is a $C(\mathscr K)$-space.
The same holds for $C_s(\mathscr K)$-spaces.
\end{proposition}

\begin{proof}
Let $F \in \mathscr K$ and  $f:F \to E \oplus L$ be nearly continuous with closed graph. From the latter follows that $f^{-1}(K)$ is closed for every compact $K \subseteq L$. Now let 
$U \subseteq L$ be open and relatively compact, then $f^{-1}(\xoverline{U})$ is closed in $F$ and contained in $f^{-1}(L)$, hence by near continuity
$f^{-1}(U) \subseteq \xoverline[.9]{f^{-1}(U)}^\circ \subseteq \xoverline[.9]{f^{-1}(U)} \subseteq f^{-1}(\xoverline{U})\subseteq f^{-1}(L)$. This shows that $f^{-1}(L)$ is open in $F$, and moreover,
since $L$ as a locally compact space is regular, that $f$ is continuous at every $x'\in f^{-1}(L)$.

 But now $f^{-1}(E)=F \setminus f^{-1}(L)$ is closed, hence again by near continuity $f^{-1}(E)$ is also open in $F$, and we have $F = f^{-1}(E) \oplus f^{-1}(L)$.  
By  hypothesis on $\mathscr K$ we have $f^{-1}(E) \in \mathscr K$,
hence $g:=f|f^{-1}(E)$ is continuous, because $G(g) = G(f) \cap (f^{-1}(E) \times E)$ is closed in $f^{-1}(E) \times E$, and the restriction $g$ is also nearly continuous as a map $g:f^{-1}(E) \to E$.
Altogether that shows continuity of $f$. In the $C_s(\mathscr K)$-case, we have just to observe that $g=f|f^{-1}(E)$ is separating when $f$ is. 
\hfill $\square$
\end{proof}

\section{$F$-heredity}
\label{sect_F}
In this section we consider invariance of $C(\mathscr K)$- and $C_s(\mathscr K)$-spaces under taking closed subspaces. Our first observation is
\begin{proposition}
\label{hereditary}
Closed subspaces of  $C(\mathscr K)$-spaces are 
$C(\mathscr K)$-spaces, and closed subspaces of $C_s(\mathscr K)$-spaces are $C_s(\mathscr K)$-spaces.
\end{proposition}

\begin{proof}
1)
Let $E$ be a $C(\mathscr K)$-space and $F \subseteq E$  closed. 
Let $f:G \to F$ be nearly continuous  with closed graph and $G \in \mathscr K$. We have to show that $f$ is continuous.
Since $F$ is closed, the graph of $f$ is also closed as a subset of $G \times E$.
Since $f$ is also nearly continuous when regarded as a mapping from $G$ to $E$, $f$ is continuous as a mapping from $G$ to $E$,
and then also as a  mapping from $G$ to $F$.

2) Now let $E$ be a $C_s(\mathscr K)$-space. Let $G \in \mathscr K$ and  $f:G \to F$ is nearly continuous and separating. 
All we have to show is that $f$ is also separating as a mapping from $G$ to $E$. 
Let $x,y\in E$, $x\not=y$. There are two cases to be considered. (i) $x,y\in F$, (ii) one of the points, say $y$, is not in $F$. In the latter case we choose
$V_y$ open such that $y\in V_y$ and $V_y \cap F =  \emptyset$, then $f^{-1}(V_y)=\emptyset$, and we are done.
In case (i) we use the fact that $f$ is separating as a mapping $G$ to $F$. The latter says that there exist open sets $V_x,V_y$ in $E$ with $x\in V_x$ and 
$y\in V_y$ such that
$\xoverline[.9]{f^{-1}(V_x \cap F)} \cap f^{-1}(V_y \cap F) = \emptyset = f^{-1}(V_x \cap F) \cap \xoverline[.9]{f^{-1}(V_y \cap F)}$. But then
since $f$ maps into $F$, this shows $x,y$ are separated.
\hfill $\square$
\end{proof}

This allows to draw  a first line between $C$- and $C_s$-spaces on one end,  and $B_r$-spaces on the other. 
We have the following:

\begin{proposition}
\label{every}
Every semi-regular space $F$ is a closed subspace of a $B_r$-space $E$ containing a dense completely metrizable subspace. When $F$ is not a
$C(\mathscr K)$-space, then  $E$ is not a $C(\mathscr K)$-space, and the same is true with regard to $C_s(\mathscr K)$-spaces.
\end{proposition} 

\begin{proof}
1)
Let
$E = (F \times \{1\}) \cup (F \times \{2\})$ and equip it with a topology as follows. For 
$U$ open in $F$ and $Y$ finite let $W(U,Y)= \{(x,2): x\in U\} \cup \{(y,1): y\in U\setminus Y\}$.
Note that $W(U_1,Y_1) \cap W(U_2,Y_2) = W(U_1 \cap U_2, Y_1 \cup Y_2)$. Now let $\tau$ be the topology on $E$ generated by the 
singletons $\{(x,1)\}$ together with the $W(U,Y)$, where $U$ is open in $F$ and $Y$ is finite.  
The space $E$ is sometimes called the Alexandrov duplicate of $F$.

2)
Let $F = G \cup H$, where $G$ is the set of isolated points in $F$, $H$ the non-isolated points. 
We show that $D:= (F \times \{1\} ) \cup (G \times \{2\}) $ is dense in $(E,\tau)$. For this it suffices to show that for every $x \in H$ any basic neighborhood 
$W(U,Y)$ of  $(x,2)$ meets $D$.  Now
$x\in U$ and since $x\in H$ is non-isolated, $U$ is infinite. Therefore $U$ contains a point $z$ with $z\not \in Y$.
Hence $(z,1) \in W(U,Y)$, and that proves the claim, because $(z,1)\in D$.
Note that $D$ is discrete in  the induced topology.

3)
Note that $F \times \{2\}$ is closed in $E$, because $F \times \{1\}$ is open. Moreover, $F \times \{2\}$ is homeomorphic to $F$,
because $W(U,Y) \cap ( F \times \{2\} )= U \times \{2\}$. 

4) We have to show that $\tau$ is a Hausdorff topology. Let $x,y\in F$, $x\not= y$. Then $(x,1)$, $(y,1)$ are clearly separated by their singletons.
Also $(x,1)$ and $(y,2)$ are separated by $\{(x,1)\}$ and $W(U_y,\{x\})$ for any neighborhood $U_y$ of $y$ in $F$. But
$(x,2)$, $(y,2)$ are also separated, because here we choose open $x \in U_x$, $y\in U_y$ in $F$ with $U_x \cap U_y=\emptyset$,
then $W(U_x,\emptyset) \cap W(U_y,\emptyset)=\emptyset$. It remains to separate $(x,1)$ and $(x,2)$. Here
we  take $\{(x,1)\}$ and $W(U_x,\{x\})$ for an open $U_x$ containing $x$.

5)
We show that $E$ is semi-regular. Since regularity at points $(x,1)$ is clear,
consider  $(x,2) \in W(U,Y)$. That means $x\in U$. Using semi-regularity of $F$ choose a regular-open $V$ in $F$ with $x \in V = \xoverline{V}^\circ \subseteq U$.
It suffices to observe that $W(V,Y)$ is regular open in $E$. 
Since $W(V,Y) = (V\times \{1\}) \cup (V \times \{2\}) \setminus (Y \times \{1\})$ and $Y$ is finite, it suffices to show that
$(V\times \{1\}) \cup (V \times \{2\})$ is regular-open, but that is clear given that $V$ is regular-open.

6) Now since $E$ contains the open dense completely metrizable subspace $D=(F \times \{1\})\cup (G \times \{2\})$ and is semi-regular, it is a $B_r$-space; cf.
\cite{open_mapping,on_the_theory}.  By Proposition \ref{hereditary}, however,   if $F$ is not a $C(\mathscr K)$-space, then $E$ cannot be a $C(\mathscr K)$-space, and similarly for the $C_s(\mathscr K)$-case.
\hfill $\square$ 
 \end{proof}
\begin{corollary}
\label{cor3}
The class of $B_r$-spaces is not F-hereditary.
\hfill $\square$
\end{corollary}

A similar counterexample can be obtained by results of Strecker and Wattel, Banaschewski,  and Berry, see \cite{strecker}, showing when combined that every
 Hausdorff space $E$ is the closed subspace of a $H$-minimal space. By Proposition \ref{hereditary} this $H$-minimal space, while a $B_r$-space,  cannot be a $C_s(\mathscr K)$-space if $E$ 
 fails to be one.

\begin{remark}
 It is interesting to compare Corollary \ref{cor3} and Proposition \ref{hereditary} to the situation in the classical categories. In Pt\'ak's theory it is known that
 $B_r$- and $B$-completeness are $F$-hereditary  (see \cite{koethe,ptak}). The same is true for $B_r$- and $B$-completeness in topological vector spaces
 (see \cite{adasch}). This can be further extended to abelian groups (see \cite[p. 46]{husain}). In general topological groups the situation is more complicated, and only
 partial results on $F$-heredity are known (see \cite{Gr}). 
 \end{remark}

\begin{corollary}
\label{cor4}
$B_r$-spaces  are not invariant under perfect preimages.
\end{corollary}

\begin{proof}
Every clopen subspace $G$ of a $B_r$-space $E$ is again $B_r$, because $E = G \oplus E\setminus G$. But if a property is inverse perfect invariant and clopen invariant, then it is $F$-hereditary
(cf. \cite[Thm. 3.7.29]{engelking}). Hence
$B_r$ cannot be inverse perfect invariant. 
\hfill $\square$   
\end{proof}

Let $\mathscr K_k$ be  $k$-spaces, 
then for the same reason the class of $B_r(\mathscr K_k)$-spaces is not inverse perfect invariant. However, by Propositions \ref{perfect} and
\ref{prop4},
 $C_s(\mathscr K_k)$-spaces are a subclass of  $B_r(\mathscr K_k)$ which {\it is} inverse perfect invariant. 
Yet another consequence of Proposition \ref{hereditary} is the following

\begin{proposition}
Let $E$ be a $C_s$-space. Then every subspace of $E$ is semi-regular.
\end{proposition}

\begin{proof}
Every closed subspace $F$ of $E$ is a $C_s$-space, hence is $B_r$, and it is known \cite{open_mapping} that every $B_r$-space is semi-regular.
Hence $F$ is semi-regular. But semi-regularity is preserved by dense subspaces \cite[2.7.6]{engelking}, hence every subspace of $E$ is semi-regular.
\hfill $\square$ 
\end{proof}

We will get back to the question of regularity of $C_s$-spaces in Section \ref{sect_regular}.

\section{Product invariance}
\label{sect_product}
Invariance of open mapping and closed graph spaces under products fails in all classical categories, and we therefore expect 
negative results in the topological context. For open mapping spaces it is known that the product of even two $B_r$-spaces need not be $B_r$, and
that there exists a $B_r$-space whose square $E \times E$ is no longer $B_r$,   cf. \cite{sums_and_products}.

For closed graph spaces $C(\mathscr K)$ and $C_s(\mathscr K)$ we have presently only scarce results. 
As seen in Section \ref{sect_perfect}, 
the product $E \times K$ of a $C(\mathscr K)$-space and a compact space is again $C(\mathscr K)$ when $\mathscr K$ is a class of $k$-spaces,  
and  $E \times D$ is $C(\mathscr K)$ if $E$ is $C(\mathscr K)$ and $D$ is discrete,
as follows from Corollary \ref{cor1}, and similarly in the $C_s(\mathscr K)$-case. But these appear to be fairly exceptional situations.

A theorem due to Mahavier, see also \cite[Thm. 5]{piotrowski} or \cite{berner}, goes as follows:

\begin{proposition}
Let $\mathscr F$ be the class of Fr\'echet spaces. Then every regular locally countably compact space $E$ is a $C(\mathscr F)$-space.
\end{proposition}

This can be used to obtain counterexamples. In the following let $\mathscr M$ and $\mathscr M_c$ stand for  metrizable, respectively,  completely metrizable spaces.

\begin{proposition}
The classes $C(\mathscr F)$, $C(\mathscr M)$ and $C(\mathscr M_c)$  are not closed under finite products.
\end{proposition}

\begin{proof}
Frol\'ik \cite{frolik_countably_compact} shows that every separable metrizable space $E$ is a closed subspace of a product of two completely regular countably compact spaces $F,G$.
Both $F,G$ are $C(\mathscr F)$-spaces by Mahavier's theorem, but their product $F \times G$  is not, as otherwise $E$ as a closed subspace would have to be $C(\mathscr F)$ 
by Proposition \ref{hereditary}, which it is not in general. Naturally, the same  argument works for the subclasses $\mathscr M$
and $\mathscr M_c$ of $\mathscr F$.
\hfill $\square$
\end{proof}

\begin{remark}
This gives also a negative answer to product invariance of  $C_s(\mathscr F)$-, $C_s(\mathscr M)$- and $C_s(\mathscr M_c)$-spaces. 
For $C$- and $C_s$-spaces product invariance is an open problem, the expected answer being in the negative. While it is folklore that locally compact spaces are $C$-spaces (cf. Proposition \ref{prop10}), 
it is not even known whether the product
$E \times L$ of a $C$-space $E$ and a locally compact space $L$ is again a $C$-space, and similarly in the $C_s$-case. 
\end{remark}

Product invariance plays a crucial role in Banach's method to derive the closed graph theorem from the open mapping theorem, which indicates that its applicability is limited.
Banach's argument can be presented under the following form.

\begin{proposition}
Let $\mathscr C$ be an $F$-hereditary class closed under taking finite products, and satisfying $\mathscr C \subseteq B_r(\mathscr C)$. Then  $\mathscr C  \subseteq
C(\mathscr C)$.
\end{proposition}

\begin{proof}
Let $E,F \in \mathscr C$ and $f:F \to E$  nearly continuous with closed graph. We have to show that $f$ is continuous. As a closed
subspace of $F \times E$, the graph $G(f)$ belongs to the class 
$\mathscr C$, hence is a $B_r(\mathscr C)$-space. Now $\phi:G(f) \to F$ defined by
$\phi(x,f(x)) = x$ is continuous since $\phi=p_F|G(f)$ with $p_F$ the projection $F \times E \to F$,
and moreover $\phi$ maps bijectively onto $F$. In addition, $\phi$ is also nearly open, because if $U$ is open in $F$ with $x\in U$ and $V$ is open
in $E$ with $f(x)\in V$, then $\phi( (U \times V) \cap G(f)) =  U \cap f^{-1}(V)$, hence $\xoverline[0.99]{\phi( (U \times V) \cap G(f)) }=
\xoverline[.99]{U \cap f^{-1}(V)} \supseteq U \cap \xoverline[.9]{f^{-1}(V)} \supseteq U \cap \xoverline[.9]{f^{-1}(V)}^\circ \supseteq U \cap f^{-1}(V)$ by near continuity of $f$.
This shows $\phi((U \times V)\cap G(f)) \subseteq \xoverline[0.99]{\phi( (U \times V) \cap G(f)) }^\circ.$

Now since $G(f)$ is a $B_r(\mathscr C)$ space and $F \in \mathscr C$, $\phi$ is open. Hence $\phi^{-1}:F \to G(f)$ is continuous, and since the projection $p_E:F \times E \to E$ 
is continuous, so is $f=p_E \circ \phi^{-1}$. 
\hfill $\square$
\end{proof}

\section{$C$- and $C(\mathscr M)$-spaces}
\label{sect_M}
The discussion so far suggests that the class of $C$-spaces is rather small, whereas $C_s$-spaces are reasonable. Proposition \ref{prop10} 
shows that locally compact spaces are $C$-spaces, but are there any other, non-locally compact
$C$-spaces?

A  set $H\subseteq E$ is separated 
if there exists a disjoint family of open sets $U_y$, indexed by $y\in H$, such that $H \cap U_y = \{y\}$. In a regular space $E$, every countable discrete $H$ is separated. 
For general discrete $H$ this holds when $E$ is collectionwise Hausdorff.
The following extends \cite[Thm. 2.3]{berner}:

\begin{proposition}
\label{isolated}
Suppose the set of non-isolated points in $E$ is separated. Then $E$ is a $C$-space.
\end{proposition}

\begin{proof}
1) Let us first consider the case where $E$ has at most two non-isolated points.
That is,  $E= G \cup \{y,z\}$, where for every $x\in G$ the singleton $\{x\}$ is open, and where $y,z$ are non-isolated. Let $f:F \to E$ be nearly continuous with closed graph.
The latter implies that all fibers $f^{-1}(x)$  are closed. By near continuity, for the $x\in G$ we have $f^{-1}(x) \subseteq \xoverline[.9]{f^{-1}(x)}^\circ = f^{-1}(x)^\circ$, hence every fiber
$f^{-1}(x)$, $x\in G$,  is closed-and-open.  In particular, $f$ is clearly continuous at every $x'\in f^{-1}(G)$.

Now let $U$ be an open neighborhood of $y$, and let $y' \in f^{-1}(y)$. 
We have to show that $f^{-1}(U)$ is a neighborhood of $y'$. By near continuity $\xoverline[.9]{f^{-1}(U)}$ is a neighborhood of $y'$,
hence there exists an open $W$ in $F$ with $y'\in W \subseteq \xoverline[.9]{f^{-1}(U)}$. 
Moreover, since $f^{-1}(z)$ is closed and $y'\not\in f^{-1}(z)$, we may assume $W \cap f^{-1}(z)=\emptyset$. We prove $W \subseteq f^{-1}(U)$.
Let $x' \in W$, then either $f(x') = y \in U$, or $f(x')  = x\in G$, in which case $f^{-1}(x)$ is a neighborhood of $x'$, which due to $x'\in \xoverline[.9]{f^{-1}(U)}$
satisfies
 $f^{-1}(x) \cap f^{-1}(U)\not=\emptyset$. That
implies $f^{-1}(x) \subseteq f^{-1}(U)$, hence $x'\in f^{-1}(U)$. 
Clearly $W \subseteq f^{-1}(U)$  proves continuity at $y'$. 
Continuity at $z'\in f^{-1}(z)$ is analogous. 

2) Concerning the case $E = G \cup H$ with $G$ the set of isolated points and $H$ the non-isolated points, by hypothesis there exist a
disjoint family $U_y$ of open sets indexed by $y\in H$ such that $U_y \cap H = \{y\}$ for every $y\in H$. Then $E \setminus \bigcup \{U_y:y\in H\} \subseteq G$
is discrete. Adding this set to one of the $U_y$
gives therefore a disjoint open cover $(U_y)$ of $E$, which means $E = \bigoplus_{y\in H} U_y$, where every space $U_y$ has precisely one non-isolated point $y$. 
The sum $U_y \oplus U_z$ of two of those has two non-isolated points, hence is a $C$-space by part 1). Then $E$ is a $C$-space by Lemma \ref{iota}. 
\hfill $\square$
\end{proof}

\begin{example}
This is useful, because even a space with just one non-isolated point need not be locally compact. We can take
the space $M_1 = \{\frac{1}{n} + \frac{1}{n^2m}: n,m\in \mathbb N\} \cup \{0\}$, then $M_1$ is a non locally compact scattered  countable
metric space with exactly one non-isolated point. It is known that {\it every} countable non locally compact scattered metric space with exactly one non-isolated point is homeomorphic to 
$M_1$; cf. \cite{gillam}.
\end{example}

\begin{remark}
If $E$ is countable with precisely one non-isolated point which does not have a countable neighborhood base, then  $E$
is a $C$-space which,
according to  \cite[Cor. 8.3]{michael_game}, is paracompact, but is neither \v Cech-complete nor sieve-complete.
\end{remark}

Note that polish spaces are $C(\mathscr B)$-spaces, $\mathscr B$ the  Baire spaces in the sense of category, and hence are
$C(\mathscr M_c)$-spaces; cf. \cite{moors1,graph_theorem} and Section \ref{sect_CM}.  Polish spaces  are also $C_s$-spaces; cf. \cite{piotrowski}, \cite[Thm. 2]{moors1} or Section \ref{sect_tandem}.  
In contrast,  for a 
polish space to be a $C(\mathscr M)$-space  we expect rather severe restrictions. 

\begin{proposition}
Let $E$ be a polish space which is a $C(\mathscr M)$-space. Then $E$ is a $K_\sigma$-space.
\end{proposition}

\begin{proof}
Assume on the contrary that  $E$ is not a $K_\sigma$, then it contains a closed copy of
the irrationals $\mathbb P$ by the Hurewicz theorem. Now as $E$ is a $C(\mathscr M)$-space, $E \oplus E$ is also $C(\mathscr M)$ 
by Proposition \ref{oplus}. Hence we have a $C(\mathscr M)$-space  $E \oplus E$ which contains a closed copy of 
$\mathbb P \oplus \mathbb P$. Suppose the point set of $E \oplus E$  is $(E \times \{1\} )\cup (E \times \{2\})$
and let $P_1$ be the closed copy of $\mathbb P$ in $E\times \{1\}$, $P_2$ the closed copy of $\mathbb P$ in $E \times \{2\}$.

Now we use a construction due to Wilhelm \cite{wilhelm_pettis}.
Choose two dense $G_\delta$ subsets $G_1,G_2$ of $\mathbb R$ such that $G_1 \cup G_2 = \mathbb R$ and $G_1\setminus G_2$,
$G_2\setminus G_1$ are both dense in $\mathbb R$, by letting $Q_1,Q_2$ be two disjoint copies of $\mathbb Q$ both dense in $\mathbb R$, and putting
$G_i = \mathbb R \setminus Q_i$. Notice that $G_1 \simeq G_2 \simeq \mathbb P \simeq P_1 \simeq P_2$.
Let $H = G_1\cap G_2$, which is also dense in $\mathbb R$, and let $H_1 \simeq H_2 \simeq H$ be two copies of $H$.

Now put $F = (G_1\setminus G_2) \cup (G_2\setminus G_1) \oplus H_1 \oplus H_2$, where all three summands have the subspace topology of $\mathbb R$, and define a mapping $f:F \to P_1 \oplus P_2 \subseteq E \oplus E$
as follows. Send $x\in G_1\setminus G_2$ to its copy $x_1\in P_1$,  $x\in G_2\setminus G_1$ to its copy $x_2\in P_2$. Send $(y,1) \in H_1$ to $y_1\in P_1$ and $(y,2)\in H_2$
to $y_2\in P_2$.  Then $f$ is bijective onto $P_1 \oplus P_2$ and as shown in \cite{wilhelm_pettis} has closed graph in $F \times (P_1 \oplus P_2)$, hence also in $F \times (E\oplus E)$,
$P_1 \oplus P_2$ being closed. One also checks that $f$ is nearly continuous, using denseness of $G_1\setminus G_2,G_2\setminus G_1$ and $H_1,H_2$ in $\mathbb R$. Since $f$ is clearly not continuous and $F \in \mathscr M$, this contradicts the fact that
$P_1 \oplus P_2$, as a closed subspace of $E \oplus E$, ought to be a $C(\mathscr M)$-space by Proposition \ref{hereditary}.
\hfill $\square$
\end{proof}

Naturally, the result could also be read in the sense that any space $E$ which contains a closed copy of $\mathbb P$ cannot be a $C(\mathscr M)$-space.

\begin{proposition}
Every polish $C(\mathscr M)$-space $E$ contains a dense open locally compact subspace $G$ and $H=E\setminus G$ is a closed nowhere dense $K_\sigma$-set of non-isolated points of $E$.
\end{proposition}

\begin{proof}
Let $E$ be a polish $C(\mathscr M)$-space, $G$ its  set of isolated points, $E = G_0 \cup H_0$, where $H_0$ is closed perfect. Suppose $H_0$ is somewhere dense in $E$. Being itself a $C(\mathscr M)$-space
by Proposition \ref{hereditary} shows $H_0$ must be $K_\sigma$, hence some compact subset $K$ of $H_0$ is somewhere dense in $E$, hence has nonempty interior. 
That means there exists a relatively compact nonempty open set $U_1 \subseteq H_0$. We put $G_1 = G_0 \cup U_1$, then $G_1$ is open in $E$ and locally compact as a subspace. 
Now we can repeat the procedure, because $H_1 = E \setminus G_1$ is closed. Hence it is either nowhere dense in $E$, or if somewhere dense, and being a $C(\mathscr M)$-space, contains again a relatively
compact open $U_2$. We put $G_2 = G_1 \cup U_2 = G \cup U_1 \cup U_2$, which is then open in $E$ and locally compact. Continuing in this way via transfinite induction
leads to a decomposition
$E = G \cup H$, where $G$ is open dense and locally compact, $H$ closed nowhere dense in $E$. Being a polish $C(\mathscr M)$-space, $H$ is a $K_\sigma$.
\hfill $\square$
\end{proof}

We cannot expect $E$ to be locally compact altogether, as the example of $M_1$ shows. But we can iterate the procedure on $H$, as this is also a polish $C(\mathscr M)$-space. 
This raises the following question: Suppose a polish space $E$ has an open dense locally compact subspace $G$, and $H=E \setminus G$ is a $C(\mathscr M)$-space in the induced topology (which is therefore a $K_\sigma$).
Does this suffice to conclude that $E$ is a $C$-space?

The following can be proved with a slight modification of the proof of Proposition \ref{isolated}.

\begin{proposition}
Suppose $E$ is regular and has a co-finite locally compact subspace $G$. Then $E$ is a $C$-space.
\end{proposition}

It is not  clear whether this continues to hold in the case where $H=E\setminus G$ is separated and infinite. 
 What about closed copies of $\mathbb Q$? 
Here the answer is clear.

\begin{proposition}
Every metrizable $C_s(\mathscr M)$-space is F-hereditary Baire. The same is true for metrizable $C(\mathscr M)$-spaces.
\end{proposition}

\begin{proof}
By a result of Morita; cf. \cite[4.4.J]{engelking},
every metrizable space $E$ is the perfect image of a  subspace $S$ of a space $D^\omega$, where $D$ is discrete. Therefore
when $E$ is $C(\mathscr M)$ or $C_s(\mathscr M)$, then so is $S$ by Proposition \ref{perfect}. In consequence $S$ is also a $B_r(\mathscr M)$-space by Proposition \ref{prop4}, 
and so is every closed subspace of $S$ by Proposition \ref{hereditary}. However, $S$ has closed-and-open basic sets, hence by \cite[Prop. 1]{baire_category} $S$
is a Baire space, and the same applies to  every closed subspace of $S$. Note that the method of proof in \cite{baire_category}  constructs for the given zero-dimensional metrizable topology $\tau$ on $S$ a coarser
topology  $\tau'$ such that $id_S:(S,\tau) \to (S,\tau')$ is nearly open, and one has but to remark that $\tau'$ is also metrizable, in order to apply the hypothesis. 
To conclude, the closed continuous image $E$ of the paracompact
$F$-hereditary Baire space $S$ is $F$-hereditary Baire, as shown in \cite[Thm. 4.10 (iii)]{haworth}, hence $E$ is $F$-hereditary Baire. 
\hfill $\square$
\end{proof}

\begin{corollary}
Let $E$ be a metrizable $C_s$-space. Then $E \oplus F$ is a $B_r$-space for every \v Cech-complete space $F$, and for every $B_r$-space $F$ analytic in the sense of Frol\'ik
{\rm \cite{frolik_analytic}}.
\end{corollary}

\begin{proof}
According to  \cite[Thm. 2, Prop. 4]{sums_and_products} a completely regular  $B_r$-space $E$ is Baire if and only if the sum
$E \oplus F$ is  a $B_r$-spaces for every  \v Cech-complete space $F$.
The same statements holds with respect to testing against all $B_r$-spaces $F$   analytic in the sense of Frol\'ik.
\hfill $\square$
\end{proof}

\section{Regularity}
\label{sect_regular}
We had already observed that $B_r$-spaces $E$ are semi-regular, while proofs of the closed graph theorem consistently  require 
regularity of $E$, as witnessed by the results in \cite{berner,wilhelm_relations,piotrowski,graph_theorem,moors1,kenderov,noll_tandem},  and also by the ones we prove here. 
Moreover, $C_s$-spaces came very close to being regular, as they are hereditary semi-regular.
This raises the question whether regularity is altogether a necessary condition. We have the following partial answer.

\begin{lemma}
\label{technical}
Suppose $E$ is a $C_s$-space. 
Let $x_0 \in E$ be a point which admits an open neighborhood basis $\mathcal V$ such that
for all $V,V' \in \mathcal V$, $V \cap V'\in \mathcal V$ and $\xoverline[.9]{V \cap V'} = \xoverline{V} \cap \widebar{V'}$.
Then $E$ is regular at $x_0$.
\end{lemma}

\begin{proof}
1)
Let $\tau$ be the topology on $E$.
Suppose there exists a $\tau$-open neighborhood $V_0$ of $x_0$ such that $\xoverline{V} \not\subset V_0$ for any $V \in \mathcal V$.
To $V \in \mathcal V$ associate $V^* = V \cup \xoverline{V} \setminus V_0$ and put $\mathcal V^* = \{V^*: V \in \mathcal V\}$. Let $\mathcal U = \{U \in \tau: x_0 \not\in U\}$
and put
$\mathcal S^* = \mathcal U \cup \mathcal V^*$. Let $\tau^*$ be the topology generated by the subbasis $\mathcal S^*$. Then 
a basis for $\tau^*$ consists of finite intersections of sets from $\mathcal S^*$. But clearly $\mathcal U$ is closed under taking finite intersections, and so is $\mathcal V^*$,
because   $V_1^* \cap V_2^* = (V_1\cap V_2)^*$ due to the technical assumption on $\mathcal V$. Hence basic sets for $\tau^*$ are the $B^* = U \cap V^*$.

2)
We show that $\tau^*$ is a Hausdorff topology. Let $x,y\in E$, $x\not=y$. When $x_0 \not\in \{x,y\}$, then we can clearly separate  by
$U_x,U_y\in \mathcal U$. It remains to separate $x \not= x_0$. Choose $U\in \mathcal U$ containing $x$ and $V\in \mathcal V$ such that
$U \cap V = \emptyset$. Then $U \cap \xoverline{V}=\emptyset$, hence $U \cap V^*=\emptyset$, which proves the claim.

3) We have to  show that $id_E:(E,\tau^*) \to (E,\tau)$ is separating. For the purpose of preparation we first prove  $cl_{\tau^*}(U) \subseteq \xoverline{U}$ for the $U \in \mathcal U$. 
Indeed, let $x\in cl_{\tau^*}(U)$, then $U \cap B^*\not=\emptyset$ for every basic $B^*$ containing $x$. Let $B^*=W \cap V^*$ one such. That means $x\in W \in \tau$ with $x_0\not\in W$
and $\emptyset \not= B^* \cap U \subseteq W \cap U$, and since $W$ is arbitrary, that shows $x\in \xoverline{U}$.

Next consider $V \in \mathcal V$. We claim that $cl_{\tau^*}(V) \subseteq \xoverline{V}$. Let $x\in cl_{\tau^*}(V)$, then
$B^* \cap V \not=\emptyset$ for every basic $B^* = W \cap V'^*$ containing $x$. Then $x\in W$ and $\emptyset \not= V \cap B^* \subseteq V \cap W$, and since
$W$ was arbitrary, we get $x\in \xoverline{V}$.

Now let us prove that $id_E$ is separating. Let $x,x'\in E$, $x\not=x'$. We have to find $U,U'\in \tau$ with $x\in U$,
$x' \in U'$, such that $cl_{\tau^*}(U) \cap U' = \emptyset = U \cap cl_{\tau^*}(U')$. When $x_0 \not\in \{x,x'\}$ simply take $U,U'\in \mathscr U$ with $U \cap U'=\emptyset$,
then  $\xoverline{U} \cap U' = \emptyset = U \cap \widebar{U'}$, hence by the above also $cl_{\tau^*}(U) \cap U' = \emptyset = U \cap cl_{\tau^*}(U')$.
Finally consider the case $x\not = x_0$. Find $U \in \mathcal U$ containing $x$ and $V \in \mathcal V$ such that $U \cap V = \emptyset$. Then
also $U \cap \xoverline{V}=\emptyset$, hence by the above $U \cap cl_{\tau^*}(V) = \emptyset$. But also $\xoverline{U} \cap V = \emptyset$,
which gives $cl_{\tau^*}(U) \cap V = \emptyset$ by the above preparation. Altogether this proves $id_E$ is separating.

4) We have to show that $id_E:(E,\tau^*) \to (E,\tau)$ is nearly continuous. Since $id_E$ is clearly continuous at points $x\not= x_0$,
it remains to prove near continuity at $x_0$. We prove this by showing $V^* \subseteq cl_{\tau^*}(V)$ for every $V \in \mathcal V$.
Since $V \subseteq cl_{\tau^*}(V)$ is clear, we have to show $\xoverline{V}\setminus V_0 \subseteq cl_{\tau^*}(V)$. Let $x\in \xoverline{V}$,
and let $B^* = U \cap V'^*$ be a basic neighborhood of $x$. When $x\in V'$, then $U \cap V'$ is a $\tau$-neighborhood of $x\in \xoverline{V}$,
hence $U \cap V' \cap V\not=\emptyset$. It remains the case where
$x\in \widebar{V'}\setminus V_0$. Then $x\in U \cap \widebar{V'} \cap \xoverline{V} = U \cap \xoverline[.9]{V' \cap V}$ using the technical hypothesis, and this implies
again $U \cap V' \cap V \not=\emptyset$. Hence in both cases $B^* \cap V\not=\emptyset$. 

5) Having shown that $id_E$ is nearly continuous and separating and $\tau^*$ is Hausdorff, the fact that $(E,\tau)$ is a $C_s$-space allows to apply the closed graph theorem and
gives continuity of $id_E$, hence $\tau \subseteq \tau^*$.
Hence there exists $V\in \mathcal V$ such that $V^* \subseteq V_0$, and that implies $\xoverline{V}\setminus V_0 \subseteq V_0$, which is absurd.
\hfill $\square$
\end{proof}

\begin{proposition}
\label{A1}
Every first countable $C_s$-space is regular.
\end{proposition}

\begin{proof}
A first countable space  satisfies the technical condition 
of the Lemma.
\hfill $\square$
\end{proof}

The Lemma applies also to spaces with linearly ordered local bases, as for instance discussed in \cite{davis}, where the author calls
such spaces lob-spaces.   See also the chain local base spaces of  \cite{howes}, or the well-based spaces of \cite{leek}.

\begin{example}
Herrlich \cite{He} constructs a first countable non-compact $H$-minimal space, which is therefore a $B_r$-space, but cannot be a $C_s$-space,
because as such would have to be regular by Proposition \ref{A1}, and this is impossible, because regular $H$-minimal spaces are compact.
\end{example}

\begin{lemma}
\label{baire}
Suppose $(E,\tau)$ is a Baire space, then the topology  $\tau^*$ constructed under the hypotheses of Lemma {\rm \ref{technical}} is also Baire.
\end{lemma}

\begin{proof}
Let $\beta^*$ be a strategy for player $\beta$ in the Banach-Mazur game on
$(E,\tau^*)$; cf. \cite{revalski}.  We have to find a strategy $\alpha^*$ for player $\alpha$ beating it.
We define a corresponding strategy $\beta$ on $(E,\tau)$. We may assume that $\beta^*$ plays with basic open sets. Then using the notation in the proof of Lemma \ref{technical}, 
let $\beta^*(\emptyset) = U_1 \cap V_1^* =: W_1^* \not=\emptyset$
be the first move of $\beta^*$.  Then $U_1 \cap V_1 \not=\emptyset$, because $V_1^* \subseteq \xoverline{V}_1$. 
Now if $x_0 \not \in U_1 \cap V_1$ we simply 
let  $W_1 = U_1 \cap V_1$.  Otherwise, when $x_0 \in U_1 \cap V_1$, we let $W_1$ be a nonempty open subset of $U_1 \cap V_1$ with $x_0 \not \in W_1$.  We  put $\beta(\emptyset) = W_1$.
This means $W_1 \in \mathcal U$.

Now let $O_1 \subseteq W_1$ be a potential response of player $\alpha$ to $W_1$ in the BM-game on $(E,\tau)$. Then $O_1 \subseteq W_1 \subseteq W_1^*$, so $O_1$ can be considered a
legitimate move of player $\alpha^*$ in response to $W_1^*$, because $x_0  \not\in O_1$, hence $O_1 \in \mathcal U$. Therefore $\beta^*(W_1^*,O_1) = W_2^*= U_2 \cap V_2^* \not=\emptyset$ is defined. Then again $U_2 \cap V_2 \not=\emptyset$, and we can
put $\beta(W_1,O_1) = W_2 := U_2 \cap V_2$. 

Next if $O_2 \subseteq W_2$ is a potential move of $\alpha$, then it is also a potential move of $\alpha^*$,
hence $\beta^*(W_1^*,O_1,W_2^*,O_2) = W_3^*$ is defined.  We have $W_3^* = U_3 \cap V_3^*$, and we put $W_3 = U_3 \cap V_3$ and let $\beta(W_1,O_1,W_2,O_2) = W_3$.

Continuing in this way defines a strategy $\beta$ for player $\beta$ on $(E,\tau)$, and since $(E,\tau)$ is Baire, there exists a strategy $\alpha$ beating it.
Let $\alpha^*$ be the strategy on $(E,\tau^*)$  which makes the same moves as $\alpha$, at least in those cases where $\alpha$ plays a set from $\mathcal U$.

Let $W_1^* \supseteq O_1 \supseteq W_2^* \supseteq O_2 \supseteq \dots$ be the play of $\alpha^*$ against $\beta^*$, then by construction there is a corresponding play
$W_1 \supseteq O_1 \supseteq W_2 \supseteq \dots$ 
of $\alpha$ against $\beta$ such that $W_i^* \supseteq W_i \supseteq O_i \supseteq W_{i+1}^* \supseteq \dots$. Hence $\bigcap_n W_n^* = \bigcap_n O_n = \bigcap_n W_n \not=\emptyset$,
because $\alpha$ wins, and hence $\alpha^*$ also wins.   
\hfill$\square$
\end{proof}

Let $\mathscr B$ be the class of Baire spaces in the sense of category.
Based on Lemma \ref{baire} we get:

\begin{corollary}
\label{cor_baire}
Let $E \in \mathscr B \cap C_s(\mathscr B)$ be first countable.  Then $E$ is regular. 
\end{corollary}

It follows in particular that
every second countable $C_s$-space is metrizable, and therefore $F$-hereditary Baire.
Similarly  every second countable $E \in \mathscr B \cap C_s(\mathscr B)$ is metrizable and hence $F$-hereditary Baire.

The construction in Lemma \ref{technical} can be modified to include the following case.

\begin{proposition}
\label{cardinal}
Let $(E,\tau)$ be a $C_s$-space, or a space in $\mathscr B \cap C_s(\mathscr B)$, and let $x_0\in E$ have local character $\kappa = \chi(x_0,E)$. Suppose that
whenever $U,V$ are open sets with $x_0\in U \subseteq V$ and $\xoverline{U} \not\subset V$, then $| \xoverline{U} \setminus V | \geq \kappa$.
Then $E$ is regular at $x_0$.
\end{proposition}

\begin{proof}
 Assuming that $E$ is not regular at $x_0$, there exists an open neighborhood $V$ of $x_0$ such that $\widebar{V'} \not\subset V$ for any neighborhood $V'$ of $x_0$.
Now let $V_\alpha$, $\alpha < \kappa$ be an open neighborhood basis of $x_0$ of sets contained in $V$. Then $|\xoverline{V}_\alpha \setminus V|\geq \kappa$
for every $\alpha$. This allows us to define a sequence $x_\alpha\in \xoverline{V}_\alpha \setminus V$, $\alpha < \kappa$ by transfinite induction such that 
$\alpha \mapsto x_\alpha$ is injective. Now define
$V_\alpha^* = V_\alpha \cup \{x_\gamma: V_\gamma \subseteq V_\alpha\}$. Let $\mathcal U = \{U\in \tau: x_0\not\in U\}$ and
let $\tau^*$ be the topology generated by the subbasis $\mathcal S^* = \mathcal U \cup \{V_\alpha^*: \alpha < \kappa\}$.

We show that $\tau^*$ is Hausdorff. When $x\not= y$ and $x_0\not\in \{x,y\}$, we separate by elements of $\mathcal U$. Let $x\not= x_0$. Choose $x\in U \in \mathcal U$ and $\alpha$
with $U \cap V_\alpha = \emptyset$. Then also $U \cap \xoverline{V}_\alpha = \emptyset$, hence
$U \cap V_\alpha^* =  \emptyset$. Namely, if $x_\gamma \in V_\alpha^*$, then $V_\gamma \subseteq V_\alpha$, hence $x_\gamma \in \xoverline{V}_\gamma \subseteq \xoverline{V}_\alpha$ by the definition of
$V_\alpha^*$. Hence $\tau^*$ is Hausdorff. 

We show that $id_E:(E,\tau^*) \to (E,\tau)$ is separating.  Note that $cl_{\tau^*}(V_{\alpha}) \subseteq \xoverline{V}_\alpha$. Indeed, let $x\in cl_{\tau^*}(V_\alpha)$
and let $x\in U \in \mathcal U$, then $U \cap V_\alpha \not=\emptyset$, hence $x\in \xoverline{V}_\alpha$. Now let $U \in \mathcal U$, then
also $cl_{\tau^*}(U) \subseteq \xoverline{U}$, because any $x\in W \in \mathcal U$ is also a $\tau^*$-neighborhood of $x$. This shows $id_E$
is separating. 

Let us now show that $id_E:(E,\tau^*) \to (E,\tau)$ is nearly continuous. Since $id_E$ is continuous at points $x\not= x_0$, it remains to prove near continuity at $x_0$.
This  will follow when we show
$V_\alpha^* \subseteq cl_{\tau^*}(V_\alpha)$. For this it suffices to show that $\{x_\gamma: V_\gamma \subseteq V_\alpha\} \subseteq cl_{\tau^*}(V_\alpha)$.
Let therefore $N$ be a basic neighborhood of $x=x_\gamma$. We have to show $N \cap V_\alpha \not=\emptyset$. Now as $\mathcal S^*$ is a subbasis, 
a basic set is of the form
$x \in N = U \cap V_{\gamma_1}^* \cap \dots \cap V_{\gamma_r}^*$ for certain $\gamma_i$.
Here those $\gamma_i$ with $x\in V_{\gamma_i}$ can be included in the $\mathcal U$-part, because then $U \cap V_{\gamma_i} \in \mathcal U$. Hence we may assume that
$x \not \in V_{\gamma_1} \cup \dots \cup V_{\gamma_r}$. In other words, for every $i$ we have
$x \in \{x_\delta: V_\delta \subseteq V_{\gamma_i}\}$. But $x=x_\gamma$, and since $\gamma \mapsto x_\gamma$ is injective,
this means $V_\gamma \subseteq V_{\gamma_i}$ for all these $i$. Hence $x_\gamma \in V_\gamma^* \subseteq V_{\gamma_1}^* \cap \dots \cap V_{\gamma_r}^*$.
In consequence we may replace $N$ by the even smaller neighborhood $N' = U \cap V_{\gamma}^*$ of $x=x_\gamma$. 

Now if $x\in V_\alpha$, then we are done, so assume $x\in \{x_\delta: V_\delta \subseteq V_\alpha\}$. Remembering that $\gamma \mapsto x_\gamma$ is injective,
this means $x_\gamma = x_\delta$ for one of the $\delta$ with $V_\delta\subseteq V_\alpha$, in other words, $V_\gamma \subseteq V_\alpha$. Then
$N' \cap V_\alpha = U \cap V_\gamma^* \cap V_\alpha \supseteq U \cap V_\gamma^* \cap V_\gamma = U \cap V_\gamma$, and since
by construction $x_\gamma \in \xoverline{V}_\gamma$, we have $U \cap V_\gamma \not=\emptyset$, which completes our argument.

Having verified all the requirements, we apply the closed graph theorem, which gives continuity of $id_E$, hence $\tau \subseteq \tau^*$.
But then there exists $\alpha$ such that $V_\alpha^* \subseteq V$. Hence $\{x_\gamma: V_\gamma \subseteq V_\alpha\} \subseteq V$. But the entire sequence $x_\gamma$ 
is constructed to lie outside $V$, and this is a contradiction.

Finally, concerning the case $E \in \mathscr B \cap C_s(\mathscr B)$,  the argument of Lemma \ref{baire} clearly carries  over 
to the present situation.\hfill $\square$
\end{proof}

This leads back to a proof for the first countable case, where $\chi(x_0,E)=\aleph_0$, because in that case all $\xoverline{V}_n \setminus V$, $n < \omega$,  are necessarily infinite,
so the mapping $n \mapsto x_n$ can be made one-to-one.

\section{Tandem Michael game}
\label{sect_tandem}
In this section we prove a	 closed graph theorem which generalizes Moors' \cite[Thm. 4]{moors1}. We start by recalling the Michael game  \cite{michael_game}.

Here players $\alpha$ and $\beta$ choose successively sets $B_1 \supseteq A_1 \supseteq B_2 \supseteq \dots$ according to the following rules. 
Player $\beta$ begins by
choosing a nonempty $B_1$. Player $\alpha$ responds by choosing a nonempty relatively open $A_1 \subseteq B_1$, that is,
$A_1 = B_1 \cap U_1\not=\emptyset$ for some open $U_1$. Then $\beta$ continues by choosing a nonempty $B_2 \subseteq A_1$, to which
$\alpha$ responds by playing $A_2 = B_2 \cap U_2 \not=\emptyset$ for some open $U_2$. Player $\alpha$ wins when $\bigcap_{n\in \mathbb N} \closure[2.5]{A}_n \not=\emptyset$,
otherwise $\beta$ wins. Player $\alpha$ wins strongly if every filter $\mathscr F$ with $A_n\in \mathscr F$ for all $n$ has a cluster point
in $\bigcap_{n\in \mathbb N} \closure[2.5]{A}_n$. We assume that players have full information of the past. Then strategies for players $\alpha$, $\beta$ are maps
where $\beta(\emptyset) = B_1$, $\alpha(B_1) = A_1$, $\beta(B_1,A_1) = B_2$, $\alpha(B_1,A_1,B_2) = A_2$, etc. 
A space is called $m$-Baire if for every strategy $\beta$ there exists a strategy $\alpha$ which wins against $\beta$, and
$m^*$-Baire, if $\alpha$ wins strongly. See \cite[Sect. 11]{noll_tandem} for this definition.

Following the idea of \cite{noll_tandem} further, we define a tandem version of the Michael game, where in order to highlight the difference,
players are called $\alpha'$, $\beta'$. Player $\beta'$ starts with $B_1\not=\emptyset$, to which
$\alpha'$ responds with a nonempty relatively open $A_1 \subseteq B_1$. Then player $\beta'$ chooses a nonempty $B_1'$, to which player
$\alpha'$ responds by a nonempty relatively open $A_1' \subseteq B_1'$. In her third move player $\beta'$ now chooses a nonempty $B_2 \subseteq A_1$,
to which player $\alpha'$ responds by a nonempty relatively open $A_2 \subseteq B_2$. Then back to the primed side, where player $\beta'$
chooses a nonempty $B_2' \subseteq A_1'$, to which player $\alpha'$ responds with a nonempty relatively open $A_2' \subseteq B_2'$.  This leads to two nested sequences 
$B_1 \supseteq A_1 \supseteq B_2 \supseteq \dots$ and $B_1' \supseteq A_1' \supseteq B_2' \supseteq \dots$, but 
arranged in the following meandering way as time proceeds:
\begin{center}
\begin{tabular}{ccccccccc}
$\beta'$ && $\alpha'$ && $\beta'$ && $\alpha'$ && $\beta'$  \\
\hline
$B_1$ & $\to$ & $A_1$               &                  & $B_2$        & $\to$               & $A_2$  &                                    &$\dots$\\
            &                  &   $\downarrow$ &                  & $\uparrow$ &                                &       $\downarrow$      & &  $\uparrow$ \\
            &                  &    $B_1'$            & $\to$ & $A_1'$         &                                &              $B_2'$              &$\to$ & $A_2'$ \\
            \hline
           & & $\beta'$ && $\alpha'$ &&$\beta'$ && $\alpha'$ 
\end{tabular}
\end{center}

\noindent
Player
$\alpha'$ wins the game if both $\bigcap_n \closure[2.5]{A}_n \not=\emptyset$ and $\bigcap_n \closure[2.4]{A_n'} \not=\emptyset$. Otherwise $\beta'$ wins.
Player $\alpha'$ wins strongly when filters $\mathscr F,\mathscr F'$  with $A_n\in \mathscr F$ for all $n$ and $A_n'\in \mathscr F'$ for all $n$ have
accumulation points in  $\bigcap_n \closure[2.5]{A}_n \not=\emptyset$, respectively,  $\bigcap_n \xoverline[.9]{A_n'} \not=\emptyset$.

\begin{definition}
We say that $E$ is $m\tau$-Baire if for every strategy $\beta'$ there exists a strategy $\alpha'$ winning against $\beta'$, and we say
$E$ is $m\tau^*$-Baire, if this $\alpha'$ wins strongly.
\end{definition}

Partition complete spaces are those in which player $\alpha$
has a strong winning strategy in the Michael game  \cite{KM}. 
This turns out equivalent to $\alpha'$ having a strong winning strategy in the tandem Michael game (see \cite[Prop. 1]{noll_tandem} for the argument). 
Hence every partition complete space is $m\tau^*$-Baire, but the $m\tau^*$-Baire class is potentially larger.

Our terminology is aligned with  \cite{noll_tandem}, where the tandem version of the Banach-Mazur game was considered. 
In that case the corresponding spaces were called $\tau$-Baire and $\tau^*$-Baire (with $\tau$ for tandem). 
Since the Michael game reproduces the Banach-Mazur game when player $\beta$
plays with open sets, $m$-Baire spaces are Baire. In the same vein, $m\tau$-Baire spaces are $\tau$-Baire, and $m\tau^*$-Baire spaces are $\tau^*$-Baire.
In a metrizable space, starred and non-starred versions of each of these games agree. 
Note that a metrizable $m$-Baire space is $F$-hereditary Baire, which indicates that this is a much stronger
property than the usual Baire category assumption; see e.g.  \cite{moors_hereditary_baire}.
We prove that every $m\tau^*$-Baire space is a $C_s$-space.

\begin{theorem}
\label{theorem1}
Let $f:F \to E$ be nearly continuous and separating. Suppose $E$ is regular and $m\tau^*$-Baire.
Then $f$ is continuous.
\end{theorem}

\begin{proof}
Let $x\in F$ and let $O$ be an open neighborhood of $f(x)$. Choose an open neighborhood $W$ of $f(x)$ with $\xoverline{W} \subseteq O$.
We have to find an open neighborhood $U$ of $x$ such that $f(U) \subseteq O$.
Suppose this fails so that $f(U) \not\subset O$ for any neighborhood $U$ of $x$. Then
$x\in \xoverline[.9]{f^{-1}(E\setminus O)} \subseteq \xoverline[.9]{f^{-1}(E\setminus \xoverline{W})}$. Since $\xoverline[.9]{f^{-1}(W)}^\circ$ is a neighborhood of $x$
by near continuity, it follows that $f^{-1}(E\setminus \xoverline{W}) \cap \xoverline[.9]{f^{-1}(W)}^\circ \not=\emptyset$. Put $V_0 = E \setminus \xoverline{W}$, $W_0=W$.
Then $V_0 \cap f(\xoverline[.9]{f^{-1}(W_0)}^\circ)\not=\emptyset$.
We define a strategy $\beta'$ for the second player in the tandem Michael game.

We start by defining $\beta'(\emptyset) = B_1 =V_0 \cap f(\xoverline[.9]{f^{-1}(W_0)}^\circ)$, which is nonempty by our preparation above.
Now let $A_1 = B_1 \cap V_1 \not=\emptyset$ be a potential move of player $\alpha'$. We have to define $\beta'(B_1,A_1)$. 
Note that $B_1 =V_1 \cap V_0 \cap f(\xoverline[.9]{f^{-1}(W_0)}^\circ)$,
hence we may assume $V_1 \subset V_0$. Then
$A_1 = f(\xoverline[.9]{f^{-1}(W_0)}^\circ) \cap V_1$. 
We choose $x_1 \in f^{-1}(V_1) \cap \xoverline[.9]{f^{-1}(W_0)}^\circ$ such that
$f(x_1)= v_1 \in A_1 = B_1 \cap V_1$, then
$\xoverline[.9]{f^{-1}(V_1)}^\circ$ is a neighborhood of $x_1$, hence 
$\xoverline[.9]{f^{-1}(V_1)}^\circ  \cap f^{-1}(W_0) \not=\emptyset$. We define
$B_1' = f\left(\xoverline[.9]{f^{-1}(V_1)}^\circ \cap f^{-1}(W_0)  \right) \not=\emptyset$ and let $\beta'(B_1,A_1) = B_1'$.
We remark that $B_1' \subseteq f(f^{-1}(W_0))=W_0$.

Let $A_1' = B_1' \cap W_1\not=\emptyset$ be a potential response of player $\alpha'$. 
We have to define $\beta'(B_1,A_1,B_1',A_1')$. 
Since $A_1' \subseteq W_0$, we can assume $W_1 \subseteq W_0$. Then $B_1' = f\left(\xoverline[.9]{f^{-1}(V_1)}^\circ\right) \cap W_1$.
We choose $y_1 \in \xoverline[.9]{f^{-1}(V_1)}^\circ\cap f^{-1}(W_1)$ such that $w_1 = f(y_1) \in B_1' \cap W_1 = A_1'$. Then $\xoverline[.9]{f^{-1}(W_1)}^\circ$ is a neighborhood
of $y_1$, hence $\xoverline[.9]{f^{-1}(W_1)}^\circ \cap f^{-1}(V_1) \not=\emptyset$.  We put
$B_2 = V_1 \cap f(\xoverline[.9]{f^{-1}(W_1)}^\circ)$, then $B_2 \subseteq A_1$, so that we may let
$B_2 = \beta'(B_1,A_1,B_1',A_1')$ be our move.

Now let $A_2 = B_2 \cap V_2 \not=\emptyset$ be a potential response of $\alpha'$ to the move $B_2$. Since $B_2 \subseteq V_1$, we may assume $V_2 \subseteq V_1$. 
Therefore $A_2 = V_2 \cap f(\xoverline[.9]{f^{-1}(W_1)}^\circ)$. We choose $x_2 \in f^{-1}(V_2) \cap \xoverline[.9]{f^{-1}(W_1)}^\circ$ such that
$v_2 = f(x_2) \in A_2 = B_2 \cap V_2$. Then $\xoverline[.9]{f^{-1}(V_2)}^\circ \cap f^{-1}(W_1) \not= \emptyset$. Hence
$\emptyset \not=B_2' = W_1 \cap f(\xoverline[.9]{f^{-1}(V_2)}^\circ) \subseteq W_1 \cap f(\xoverline[.9]{f^{-1}(V_1)}^\circ) = A_1'$, so we may define
$\beta'(B_1,A_1,B_1',A_1',B_2,A_2) = B_2'$.

Continuing in this way defines a strategy $\beta'$, and since $E$ is $m\tau^*$-Baire, there exists a strategy $\alpha'$ beating $\beta'$ strongly.
Let $B_1,A_1,B_1',A_1',B_2,A_2,\dots$ be the play between $\alpha'$ and $\beta'$. By construction  there  exist
open sets $V_n,W_n$ with
\begin{itemize}
\item[i.] $B_{n} = V_{n-1} \cap f(\xoverline[.9]{f^{-1}(W_{n-1})}^\circ) \not=\emptyset$,  $A_n = B_n \cap V_n$,
\item[ii.] $B_n'= W_{n-1} \cap f(\xoverline[.9]{f^{-1}(V_{n})}^\circ)\not=\emptyset$, $A_n' = B_n' \cap W_n$,
\item[iii.] $V_{n+1} \subseteq V_n$, $V_1 \subseteq V_0 = E\setminus \xoverline{W}$, $W_{n+1} \subseteq W_n$, $W_1 \subseteq W_0=W$. 
\end{itemize} 
Since $\alpha'$ wins strongly, 
$K=\bigcap_n \closure[2.5]{A}_n \not=\emptyset$, $K'=\bigcap_n \xoverline[.9]{A_n'}\not=\emptyset$, and
$(V_n)$ is a neighborhood basis of the compact $K$, and $(W_n)$ a neighborhood base of the compact  $K'$.
Since $K \subseteq V_0=E\setminus \xoverline{W}$ and $K'\subseteq W$, we have $K \cap K' =\emptyset$. Hence
by  \cite[Lemma 3]{piotrowski} there exist neighborhoods $\widetilde{V} \supseteq K$ and $\widetilde{W} \supseteq K'$
which are separating, i.e., $f^{-1}(\widetilde{V}) \cap \xoverline[.9]{f^{-1}(\widetilde{W})}=\emptyset$. But $V_n \subseteq \widetilde{V}$ and $W_n \subseteq \widetilde{W}$ for some large enough $n$,
hence
$B_{n+1}=f^{-1}(V_n) \cap \xoverline[.9]{f^{-1}(W_n)} = \emptyset$, and that contradicts item i., because all $B_n$ are by construction nonempty. 
\hfill $\square$
\end{proof}

This extends Moors' \cite[Thm. 4]{moors1}.  It is not sufficient to assume that $E$ is regular and contains a dense completely metrizable subspace, as can be seen
from Propositions \ref{hereditary} and \ref{every}.

\section{$C(\mathscr M_c)$-spaces via the Michael game}
\label{sect_CM}
We apply the tandem Michael game again to prove that $m\tau^*$-Baire spaces are $C(\mathscr M_c)$-spaces. We can even go a little beyond $\mathscr M_c$.

A pair $(T,\phi)$ consisting of a tree
$T=(T,\leqslant_T)$  of height $\omega$ and a mapping $\phi$  from $T$ to the nonempty open sets in the space $F$ is called a {\it web} on $F$ (cf. \cite[Sect. 2]{noll_tandem})
if 
\begin{itemize}
\item[$(w_1)$] $\{\phi(t): t \in T\}$ is a pseudo-base of $F$, i.e., every non-empty open $U \subseteq F$ contains some
$\phi(t)$.
\item[$(w_2)$] For every $t\in T$ the set $\{\phi(s): t <_T s\}$ is a pseudo-base of $\phi(t)$, i.e., every non-empty open
$V \subseteq \phi(t)$ contains some   $\phi(s)$ with $t <_T s$.
\end{itemize}
The web $(T,\phi)$ is {\it $p$-complete} if for every cofinal branch $b\subseteq T$ the intersection $\bigcap_{t\in b} \phi(t) \not=\emptyset$ is non-empty.
The name is chosen in allusion to the pseudo-complete spaces in the sense of Oxtoby.
Spaces admitting a $p$-complete web are called  $p$-complete. This coincides with those spaces where
player $\alpha$ has a winning strategy in the Banach-Mazur game, \cite{white}.  Every regular pseudo-complete space in the sense of
Oxtoby is  $p$-complete.

\begin{theorem}
\label{p_complete}
Let $f$ be a nearly continuous closed graph mapping from a $p$-complete space $F$ to a regular $m\tau^*$-Baire space $E$. Then $f$ is continuous.
\end{theorem}

\begin{proof}
Let $x\in F$ and let $O$ be an open neighborhood of $f(x)$. Choose an open neighborhood $W$ of $f(x)$ with $\xoverline{W} \subseteq O$.
We have to find an open neighborhood $U$ of $x$ such that $f(U) \subseteq O$.
Suppose this fails so that $f(U) \not\subset O$ for every neighborhood $U$ of $x$. Then
$x\in \xoverline[.9]{f^{-1}(E\setminus O)} \subseteq \xoverline[.9]{f^{-1}(E\setminus \xoverline{W})}$. Since $\xoverline[.9]{f^{-1}(W)}^\circ$ is a neighborhood of $x$
by near continuity, it follows that $f^{-1}(E\setminus \xoverline{W}) \cap \xoverline[.9]{f^{-1}(W)}^\circ \not=\emptyset$. Put $V_0 = E \setminus \xoverline{W}$, $W_0=W$.
Then $V_0 \cap f(\xoverline[.9]{f^{-1}(W_0)}^\circ)\not=\emptyset$.

We define a strategy $\beta'$ for the second player in the tandem Michael game, making use of a  $p$-complete web $(T,\phi)$ on $F$.
We start by defining $\beta'(\emptyset)$. From the above we know $f^{-1}(V_0) \cap \xoverline[.9]{f^{-1}(W_0)}^\circ \not=\emptyset$,
hence also $\xoverline[.9]{f^{-1}(V_0)}^\circ \cap \xoverline[.9]{f^{-1}(W_0)}^\circ \not=\emptyset$.  Using $(w_1)$ choose
$t_0 \in T$ such that $\phi(t_0) \subseteq \xoverline[.9]{f^{-1}(W_0)}^\circ$ and $\phi(t_0) \cap f^{-1}(V_0)\not=\emptyset$. We put
$B_1 = V_0 \cap f(\phi(t_0))\not=\emptyset$ and let this  be our move $\beta'(\emptyset) = B_1$. 

Let $A_1 = B_1 \cap V_1 \not=\emptyset$ be a potential response of player $\alpha'$. We have to define
$\beta'(B_1,A_1)$. Since $B_1 \subseteq V_0$,  we may assume $V_1 \subseteq V_0$, hence $A_1 = V_1 \cap f(\phi(t_0))$. 
Then $\xoverline[.9]{f^{-1}(V_1)}^\circ \cap \phi(t_0)\not=\emptyset$. Using $(w_2)$ we choose $t_0 <_Tt_1$ with $\phi(t_1) \subseteq \xoverline[.9]{f^{-1}(V_1)}^\circ$ and
$\phi(t_1) \cap f^{-1}(W_0) \not=\emptyset$, the latter because $\phi(t_0) \subseteq \xoverline[.9]{f^{-1}(W_0)}$. 
We put $B_1' = f(\phi(t_1)) \cap W_0$
and let $\beta'(B_1,A_1) = B_1'$ be our move.

Let $A_1' = B_1' \cap W_1$ be a potential response of $\alpha'$. We have to define $\beta'(B_1,A_1,B_1',A_1')$. Since $A_1' \subseteq W_0$, we may assume $W_1 \subseteq W_0$. 
Then $A_1' = f(\phi(t_1)) \cap W_1$. It follows that $\phi(t_1) \cap \xoverline[.9]{f^{-1}(W_1)}^\circ \not=\emptyset$.  Hence using the pseudo-base property $(w_2)$ of $(T,\phi)$ again we can find
 $t_2 \in T$ with $t_1 <_Tt_2$ such that $\phi(t_2) \subseteq \xoverline[.9]{f^{-1}(W_1)}^\circ$ and
$\phi(t_2) \cap f^{-1}(V_1) \not=\emptyset$, the latter since $\phi(t_1) \subseteq \xoverline[.9]{f^{-1}(V_1)}$. 
Now
$B_2 := f(\phi(t_2)) \cap V_1 \subseteq A_1$ is non-empty, and we  let this be our move $\beta'(B_1,A_1,B_1',A_1') = B_2$. 

Next let $A_2 = B_2 \cap V_2$ be a potential response of $\alpha'$, $V_2 \subseteq V_1$. 
Then $\xoverline[.9]{f^{-1}(V_2)}^\circ \cap \phi(t_2)\not=\emptyset$.  Choose $t_3$ with $t_2 <_Tt_3$ such that $\phi(t_3) \subseteq \xoverline[.9]{f^{-1}(V_2)}^\circ$ and
$\phi(t_3) \cap f^{-1}(W_1)\not=\emptyset$. 
Putting $B_2' = f(\phi(t_3)) \cap W_1$ and noting that $B_2' \subseteq A_1'$, we let our move be $\beta'(B_1,A_1,B_1',A_1',B_2,A_2)=B_2'$.

Continuing in this way defines a strategy $\beta'$, and since $E$ is $\beta'$-defavorable in the strong sense in the tandem Michael game, there exists a strategy $\alpha'$
beating $\beta'$ strongly. Let $B_1,A_1,B_1',A_1',B_2,A_2,\dots$ be their play, then by construction there exist nested sequences $(V_n)$, $(W_n)$ of open sets in $E$
and a cofinal branch $t_0 <_T t_1 <_T t_2 <_T\dots$ such that
\begin{enumerate}
\item $A_k = B_k\cap V_k$, $A_k' = B_k' \cap W_k$, $V_0 = E\setminus \xoverline{W}$, $W_0=W$, $W_k \subseteq W_{k-1}$, $V_k \subseteq V_{k-1}$.
\item $B_k = V_{k-1} \cap f(\phi(t_{2k-2})) \not=\emptyset$, $\phi(t_{2k}) \subseteq \phi(t_{2k-1}) \cap \xoverline[.9]{f^{-1}(W_k)}^\circ$.
\item $B_k' = W_{k-1} \cap f(\phi(t_{2k-1}))\not=\emptyset$, $\phi(t_{2k-1}) \subseteq \phi(t_{2k-2}) \cap \xoverline[.9]{f^{-1}(V_k)}^\circ.$
\end{enumerate} 
Since $(T,\phi)$ is $p$-complete, there exists $x\in \bigcap_n \phi(t_n)$.
Now let $\mathcal N$ be the set of all pairs $(N,k)$, where $N$ is a neighborhood of $x$ contained in $\phi(t_{2k})$. By point 2. above
there exists $w(N,k)\in N \cap f^{-1}(W_k)$, and by point 3. there exists $v(N,k) \in N \cap f^{-1}(V_k)$. Let $\mathcal N$ be ordered by
$(N,k) \preceq (N',k')$ iff $N' \subseteq N$ and $k'\geq k$, then the order is directed. Since $\alpha'$ is strongly winning, the net $\langle f(v(N,k)): (N,k) \in \mathcal N\rangle$ 
has an accumulation point $v\in \bigcap_n V_n \subseteq V_0=E \setminus \xoverline{W}$, 
and the net $\langle f(w(N,k)): (N,k)\in \mathcal N\rangle$ has an accumulation point $w \in \bigcap_n W_n \subseteq W_0=W$. On the other hand, both nets $\langle v(N,k):(N,k)\in \mathcal N\rangle$ and $\langle w(N,k):(N,k)\in \mathcal N\rangle$ converge to
$x$. Since the graph of $f$ is closed, we deduce $(x,v)\in G(f)$ and $(x,w)\in G(f)$, and that implies $v=w$, hence $v=w\in V_0 \cap W_0$, which is absurd.
This completes the proof.
\hfill $\square$
\end{proof}

One deduces directly that
every nearly continuous closed graph mapping $f$ from a space $F$ containing a dense completely metrizable space into a regular $m\tau^*$-Baire space $E$ is continuous.
This is yet another extension of Pettis' classical closed graph theorem.

\begin{remark}
Even with $F$ completely metrizable,
the result  fails already when we only assume that $E$ is regular and contains a dense completely metrizable subspace, as
can be seen through Propositions \ref{hereditary} and  \ref{every}.
This weaker hypothesis {\it does} suffice if  $f$ is, in addition, assumed nearly open, see \cite{noll_tandem}.
\end{remark}

\begin{example}
Letting $X,Y$ two disjoint dense subsets of $\mathbb R$ with $X \cup Y = \mathbb R$ which are both Baire spaces in the induced topology,
 the identity from $F=\mathbb  R$ to $E=X \oplus Y$ is nearly continuous closed graph and not continuous. This shows that mere Baire spaces $E$ are not $C(\mathscr M_c)$.
One also checks that the identity is separating, hence Baire spaces are not $C_s(\mathscr M_c)$ either.
It is not known whether $m\tau^*$-Baire spaces are $C(\mathscr B)$-spaces, but the expected answer is in the negative.
\end{example}

\begin{example}
Consider a space $E$ with a single non-isolated point $y$. Then player $\alpha$ has a strong winning strategy in the Michael game. If the first move $B_1$ of player $\beta$
contains an isolated point $x$, then $\alpha$ plays $A_1 = \{x\}$, while if $B_1 = \{y\}$ for the non-isolated point $y$, then $\beta$ looses anyway. However, the situation is different for the Choquet game, 
(cf. \cite{choquet,debs}), where
player $\beta^c$ plays an open set $U_1$ and a point $x_1\in U_1$, to which player $\alpha^c$ has to respond with an open set
$V_1$ with $x_1 \in V_1 \subseteq U_1$, etc. Here player $\alpha^c$ has a winning strategy, but not necessarily a strategy winning strongly. Namely, if  
the isolated point does not have a countable neighborhood basis, then a filter $\mathcal F$ containing all $V_k$ will not necessarily have the non-isolated point
as accumulation point; cf. 
\cite[Cor. 8.3]{michael_game}.
\end{example}

\section{$C(\mathscr M_c)$-spaces}
\label{sect_final}
Kenderov and Revalski \cite{kenderov} 
prove that a regular space $E$ contains a dense completely metrizable subspace if and only if
it satisfies a variant of the closed graph theorem for  densely defined set-valued mappings $f:F \rightrightarrows E$ from Baire spaces
$F$ into $E$, where for the necessary part it suffices to test against completely metrizable $F$.  The role of near continuity of a function $f$ is played by a 
suitable multi-valued version called lower demi-continuity.
Variants of the closed graph theorem which,  in the spirit of \cite{kenderov},  derive continuity on a dense set of points, or  feeble continuity 
from near feeble  continuity,
are also known in the function case. See for instance \cite{cao,preservation,completeness_concept,KMR}. 

In the following
we take \cite{kenderov} as a motivation to obtain an
internal characterization of $C(\mathscr M_c)$-spaces.
We call $(T,\phi)$ a sieve if $T$ is a tree of height $\omega$ and $\phi$ maps $T$ to the nonempty subsets of $E$
such that the following are satisfied:
\begin{align*}
(s_1) &  \hspace{.2cm}\mbox{$E= \bigcup\{\phi(t): t\in T\}$.}\\
(s_2) &  \hspace{.2cm}\mbox{For every $t\in T$, $\phi(t) = \bigcup\{\phi(s): t <_T s \in T\}$.} \hspace{4.99cm}
\end{align*}
The sieve is $p$-complete if for every cofinal branch $b \subseteq T$, $\bigcap_{t\in b} \phi(t) \not=\emptyset$.
It is a $\delta$-sieve if $P(b):= \bigcap_{t\in b} \xoverline{\phi(t)}$ contains at most one point. In a $p$-complete $\delta$-sieve, therefore,
$\bigcap_{t\in b} \phi(t) = \bigcap_{t\in b} \xoverline{\phi(t)} = \{y_b\}$ for a single point.

\begin{lemma}
Let $f$ be a continuous surjection from a complete metric space $F$ onto a regular space $E$. Then $E$ has a $p$-complete $\delta$-sieve with the following property:
\begin{align*}
(\mu) \quad 
\begin{array}{l}\mbox{\rm For every $y\in  U \subseteq E$ open and $b\subseteq T$ a cofinal branch}\\
                     \mbox{\rm with $y\in P(b)$  there exists $t_0\in b$ such that $\phi(t_0) \subseteq U$} 
                     \end{array}
\end{align*}
\end{lemma}

\begin{proof}
Let $T$ be the tree of finite sequences of nonempty open sets in
$F$ such that $V_1 \supseteq \xoverline{V}_2 \supseteq V_2 \supseteq \xoverline{V}_3 \supseteq \dots \supseteq V_k$ and diam$(V_k) \leq 1/k$ with respect to some fixed complete metric
on $F$, ordered by extension of sequences, and define $\phi:(V_1,\dots,V_k) \mapsto f(V_k)$. Clearly for every cofinal branch, represented by a nested sequence $(V_k)$, we have
$\bigcap_k V_k = \{x\}$ for some $x\in F$, hence $y=f(x) \in \bigcap_k f(V_k)$. Suppose $\bigcap_k \xoverline[.9]{f(V_k)}$ contains another  point $y'\not= f(x)=y$, then
$(x,y') \not\in G(f)$, hence $G(f)$ being closed, there exists an open ball $B(x,\eta)$ and a neighborhood $U'$ of $y'$ in $E$ such that $(B(x,\eta) \times U') \cap G(f)=\emptyset$. But 
$V_k \subseteq B(x,\eta)$ for  $k$ large enough, hence $f(V_k) \cap U' = \emptyset$, which shows $y' \not\in \xoverline[.9]{f(V_k)}$. 

The property  $(\mu)$ can be seen as follows. Let $\{y\} = \bigcap_k f(V_k) = \bigcap_k \xoverline[.9]{f(V_k)}$ for the cofinal branch given by $(V_k)$, and suppose 
$y \in U$. Find the $x$ with $\{x\} = \bigcap_k V_k$, then $y=f(x)$.
Choose $B(x,\eta)$ such that $f(B(x,\eta)) \subseteq U$ using continuity. Since $V_k \subseteq B(x,\eta)$ for $k$ large enough, it follows that $f(V_k) \subseteq U$.
\hfill $\square$
\end{proof}

We call a sieve satisfying $(\mu)$ a $\mu$-sieve.
We now modify the definition of a $\mu$-sieve as follows, calling what is obtained a quasi-$\mu$-sieve:
\begin{align*}
\mbox{$(q\mu)$}\quad 
\begin{array}{l}
\mbox{For every $y\in U \subseteq E$ open and $b\subseteq T$ a cofinal branch with}\\
\mbox{$y\in P(b)$ there exists $t_0\in b$ such that $U \cap \phi(t)\not=\emptyset$ for all $t_0 <_Tt$}
\end{array}
\end{align*}
We are now ready to state

\begin{theorem}
\label{kenderov}
The following  are equivalent:
\begin{enumerate}
\item[\rm (1)] $E$ is a $C(\mathscr M_c)$-space.
\item[\rm (2)] Every $p$-complete $\delta$-sieve $(T,\phi)$ which is a quasi-$\mu$-sieve is a $\mu$-sieve.
\end{enumerate}
\end{theorem}

\begin{proof}
1) Assume condition (2). Let $f:F \to E$ be nearly continuous with closed graph and $F\in \mathscr M_c$. We may without loss assume
that $f$ is surjective,  because if it is not, then we let $F' = F \oplus E_d$, where $E_d$ is the set $E$ with the discrete topology, and we define $f':F' \to E$
by $f' = f \oplus id_E$, which is still nearly continuous, has closed graph, and in addition, is surjective. When we can conclude that $f'$ is continuous, then so is $f$.

Now define the sieve $(T,\phi)$ as in the proof of the Lemma above, then $(T,\phi)$ is a $p$-complete $\delta$-sieve, because in that part of the proof we only used closedness of the graph. 
It remains to check that $(T,\phi)$ is quasi-$\mu$. Let $(V_k)$ be a cofinal branch with $\{y\} = \bigcap_k f(V_k) = \bigcap_k \xoverline[.9]{f(V_k)}$ and let $y \in U$. Find
the point $x$ with $\{x\} = \bigcap_k V_k$, then $y=f(x)$.
By near continuity 
$\xoverline[.9]{f^{-1}(U)}$ is a neighborhood of $x$, hence there exists $k$ such that $V_k \subseteq \xoverline[.9]{f^{-1}(U)}$. But that means $V' \cap f^{-1}(U)\not=\emptyset$ for every
potential successor $\emptyset \not=V' \subseteq \xoverline[.9]{V'} \subseteq V_k$ in the buildup of the tree. Then $f(V') \cap U \not=\emptyset$, and that is the quasi-$\mu$ condition.

By condition (2) the sieve is even a $\mu$-sieve. 
Now let $x_n$ be a sequence in $F$ with $x_n \to x$. Let $b$ be a cofinal branch, represented by $(V_k)$, for which $\bigcap_k V_k = \{x\}$, and hence $\{f(x)\} = \bigcap_k f(V_k)
= \bigcap_k \xoverline[.9]{f(V_k)}$.  For every $k$ there exists $n(k)$ such that $x_n \in V_k$ for all $n > n(k)$, hence $f(x_n) \in f(V_k)$ for all $n > n(k)$. Since we have a $\mu$-sieve,
this implies $f(x_n) \to f(x)$.

2) Let us now assume that $E$ is a $C(\mathscr M_c)$-space. Let $(T,\phi)$ be a $p$-complete quasi-$\mu$ $\delta$-sieve on $E$. We have to show that $(T,\phi)$
is even a $\mu$-sieve.

Let $\Sigma$ be the set of all cofinal branches of $T$. We define a metric $d$ on $\Sigma$ by
$$
d(b,b') = \left\{
\begin{array}{ll}
\frac{1}{n+1} & \mbox{if $b_0=b_0', \dots,b_{n-1}=b_{n-1}', b_n \not= b_n'$} \\
0 & \mbox{else}
\end{array}
\right.
$$
Then $(\Sigma,d)$ is a complete metric space due to $p$-completeness of the sieve. Let $f:\Sigma \to E$
map the cofinal branch $b\in \Sigma$ to  the unique point $x_b\in E$ with $\{x_b\} = \bigcap_{t\in b} \xoverline[.9]{\phi(t)}$. We will show that
$f$ has closed graph and is nearly continuous.

Suppose $(b,y)\not\in G(f)$. Then there exists $t_0\in b$ such that $y\not \in \xoverline[.9]{\phi(t_0)}$. Choose a neighborhood $U$ of $y$ with $U \cap \phi(t_0)=\emptyset$. 
Suppose $t_0$ is  the $n$th element of the branch $b$. 
Then by construction of the metric, $(B(b,\frac{1}{n+1}) \times U) \cap G(f)=\emptyset$. That shows closedness of the graph.

Now let us check near continuity of $f$. Let $b\in \Sigma$, $f(b)=x$, $U$ a neighborhood of $x$. By the quasi-$\mu$ property there exists $t\in b$
such that $\phi(s) \cap U \not=\emptyset$ for all $t<_Ts$. Suppose  $t$ is the $n$th element of $b$. We shall prove
$B(b,\frac{1}{n+1}) \subseteq \xoverline[.9]{f^{-1}(U)}$. Let $b' \in B(b,\frac{1}{n+1})$, then $b,b'$ agree on the places $0,\dots,n$. Let $m\in \omega$. We have to show 
$B(b',\frac{1}{m+1}) \cap f^{-1}(U)\not=\emptyset$. We may assume $m > n$. But $(b_0',\dots,b_m') = (b_0,\dots,b_n,s_1,\dots,s_k)$ and $\phi(s_k) \cap U \not=\emptyset$
by the quasi-$\mu$-property, because $b_n=t <_T s_1 <_T \dots <_T s_k$. Choose $z$ in this intersection and successors $s_k <_Tr_{k+1}<_Tr_{r+2}<_T\dots$ such that $z\in \phi(r_{k+i})$ for all $i$
using property $(s_2)$ of the sieve $(T,\phi)$. Let $b''$ be the cofinal branch
$b''=(b_0',\dots,b_m',r_{k+1},r_{k+2},\dots)$, then $f(b'') = z$ and $d(b',b'') < \frac{1}{m+1}$, and that proves our claim.

 Using the fact that $E$ is a $C(\mathscr M_c)$-space, the mapping $f$ is continuous. This implies that the sieve $(T,\phi)$ has the $\mu$-property. Indeed, let 
 $(x_n)$ be a sequence with $x_n\in \phi(t_n)$ for a cofinal branch $b=(t_n)$, then we have to show $x_n \to x$, where $\{x\} = \bigcap_k \phi(t_k)$. Let $b^n$ be a cofinal branch which 
 agrees with $b$ in the first $n$ places and otherwise fixes $x_n \in P(b^n)$. Then $b^n \to b$, and hence by continuity $x_n \to x$. 
 \hfill $\square$
 \end{proof}

\begin{remark}
Theorem \ref{p_complete} shows  that in a regular $m\tau^*$-Baire space condition (2) in Theorem \ref{kenderov} is satisfied. 
\end{remark}

\end{document}